\documentclass[review,numbers,3p,12pt,times]{elsarticle}

\usepackage{amssymb,bm,amsmath,amsfonts}
\usepackage{verbatim}
\usepackage{epsfig,bm,color}
\usepackage[nodots]{numcompress}
\usepackage[]{units}
\usepackage{longtable}
\usepackage{multirow}
\usepackage{booktabs}

\usepackage[list]{caption,subcaption}
\usepackage{caption, subcaption}
\usepackage{listings} 

\usepackage{chemist}         
\usepackage{chemformula}

\usepackage[linesnumbered,ruled,vlined]{algorithm2e}

\usepackage{epsfig,graphicx}
\usepackage{xcolor,import}
\usepackage{transparent}

\usepackage{ulem}   
\usepackage{cancel}

\usepackage[list]{caption,subcaption}
\usepackage{cleveref}
\crefname{figure}{figure}{figures}
\crefname{question}{question}{questions}

\definecolor{dodgerblue}{rgb}{0.12, 0.56, 1.0}

\SetCommentSty{mycommfont}

\SetKwInput{KwInput}{Input}                
\SetKwInput{KwOutput}{Output}              

\definecolor{commentColor}{rgb}{0,0.6,0}
\definecolor{numberColor}{rgb}{0.5,0.5,0.5}
\definecolor{stringColor}{rgb}{0.58,0,0.82}
\definecolor{bgColor}{rgb}{0.95,0.95,0.92}
\definecolor{titleColor}{rgb}{0.00,0.00,0.80}
\definecolor{underBrace}{RGB}{30,199,166}
\definecolor{darkblue}{rgb}{0.0, 0.0, 0.55}

\lstdefinestyle{CStyle}{
  backgroundcolor=\color{bgColor},
  belowcaptionskip=1\baselineskip,
  numberstyle=\tiny\color{numberColor},
  stringstyle=\color{stringColor},
  keywordstyle=\color{magenta},
  commentstyle=\color{commentColor},
  basicstyle=\footnotesize\ttfamily,
  identifierstyle=\color{blue},
  breakatwhitespace=false,
  breaklines=true,
  captionpos=b,
  keepspaces=true,
  numbers=left,
  numbersep=5pt,
  showspaces=false,
  showstringspaces=false,
  showtabs=false,
  tabsize=2,
  language=C
}

\newcommand\redsout{\bgroup\markoverwith{\textcolor{red}{\rule[0.5ex]{3pt}{0.5pt}}}\ULon} 

\biboptions{sort&compress,square,numbers}

\journal{Mathematics in Science and Industry}
\date{\today}

\begin{document}

\begin{frontmatter}

  \title{Moving Mesh Simulations of Pitting Corrosion}

  \author[ref1]{Abu Naser Sarker }

  \cortext[cor1]{Corresponding author}
  \author[ref2]{Ronald D. Haynes \corref{cor1}}
  \ead{rhaynes@mun.ca}

  \author[ref3]{Michael Robertson}

  \address[ref1]{Scientific Computing Program, Memorial University, St. John's, NL, Canada, A1C 5S7.}
  \address[ref2]{Department of Mathematics and Statistics, Memorial University, NL, Canada, A1C 5S7.}
  \address[ref3]{Department of Physics, Acadia University, Wolfville, NS, Canada, B4P 2R6.}

  \begin{abstract}
    Damages due to pitting corrosion of metals cost industry billions of dollars
    per year and can put human lives at risk.
    The design and implementation of an adaptive moving mesh method is provided
    for a moving boundary problem related to pitting corrosion. The adaptive
    mesh is generated automatically by solving a mesh PDE coupled to
    the nonlinear potential problem. The moving mesh approach is shown to enable
    initial mesh generation, provide mesh recovery and is able to smoothly
    tackle changing pit geometry. Materials with varying crystallography are
    considered. Changing mesh topology due to the merging of pits is also
    handled. The evolution of the pit shape, the pit depth and the pit width are computed and compared to existing results in the literature.
  \end{abstract}

  \begin{keyword}
    Pitting Corrosion, Adaptive Moving Mesh, MMPDE, FEM, Crystallography
  \end{keyword}

\end{frontmatter}

{\bf Declaration of Competing Interests:} We have no competing interests to declare.

\section{Introduction}
Corrosion is a deterioration or breakdown of a material due to chemical or
electrochemical reactions.  In particular, 
pitting corrosion is one of the most disastrous and devastating localized forms of corrosion;
generating a small pit, cavity or hole in the metal. Pitting corrosion is difficult to identify, and can have a big impact on the structural integrity of metal~\cite{yu2016micromechanics,simon1999influence}.
The pit geometry depends on many factors such as the components of the metal, the surface orientation, and the physical and chemical environment at the time of attack~\cite{sharland1987review}.  Corrosion pits can have different shapes
~\cite{makhlouf2018corrosion} and with the ability to grow over time, failure of
engineering structures such as bridges, pipelines and nuclear power plants  may  result~\cite{makhlouf2018corrosion,roberge2008corrosion,cattant2008corrosion}.

Computational modeling and simulations have been a tremendous asset in the study of pitting corrosion over a wide range of conditions and materials.
Determining the pitting behavior experimentally is time        consuming,
expensive and physically difficult or           impossible in many situations.
 Numerical simulations allow us to study pitting  under a wide range  of conditions within a reasonable time.
In last few years, there have been several review papers that have focused on
partial differential equation (PDE) based models for pitting corrosion based on finite element or finite volume
methods
\cite{degiorgi2013numerical,sharland1988model2,scheiner2009finite,laycock2014computer,krouse2014modelling,de2017finite,walton1990mathematical,sharland1988model1}.
In 2019, an extensive overview of the mathematical models for pitting corrosion
based on the anodic reaction at the corrosion front, the transportation of ions
in the pits of the electrolyte domain and the pit growth over time is provided
in \cite{jafarzadeh2019computational}.  In many of these studies
COMSOL$^\text{\textregistered}$ is used to solve the PDE in the electrolyte
domain and the corrosion front movement and meshing is computed by the arbitrary
Lagrangian-Eulerian (ALE) approach and the level set
method~\cite{degiorgi2013numerical,kota2013microstructure}. In other studies, a
2D PDE model is solved with the finite element
method~\cite{duddu2014numerical,turnbull2009challenges} and the finite volume
method~\cite{scheiner2007stable,scheiner2009finite}. Pit growth is determined
by finite element methods and a level set approach
in~\cite{turnbull2009challenges}, and using an extended finite element method
(XFEM) and level set method in~\cite{duddu2014numerical}. In 2020, an ALE
method is implemented to move the mesh at the pit boundary and analyze the
relationship between the corrosion behavior and the local corrosive environment
within a single pit~\cite{wang2019multi}.

The previously mentioned FEM approaches relied on a complete remeshing of the
domain at every time step.  An alternative technique, presented here, uses an adaptive moving mesh method where the mesh size, shape and orientation of the mesh elements are automatically and continuously varied for each time step, while keeping the number of nodes and mesh topology fixed throughout the computation.

%

Continuous mesh movement approaches are divided into two main categories:
velocity-based approaches and location-based approaches. Most
velocity-based approaches are motivated by the Lagrangian algorithm, where the
mesh movement is tightly associated with the fluid or material
particle flow.  The Eulerian approach has a fixed computational mesh
and the continuum moves with respect to the mesh nodes. The  Eulerian and
Lagrangian algorithms  are commonly used in fluid
dynamics and structural material problems, respectively \cite{donea2004arbitrary}.
In general, Eulerian meshes avoid mesh tangling and diffusive solutions, but the
method can have difficulty adjusting to sharp material interfaces. One of the
advantages of the Lagrangian approach is that the advective terms do not appear
in the governing equations. Thus, the Lagrangian methods are less diffusive
compared to the Eulerian approach, while also maintaining sharp material
interfaces~\cite{huang2010adaptive}. The ALE
methods are velocity-based methods, which provide a combination of Lagrangian and Eulerian approaches~\cite{hirt1974arbitrary,hirt1997arbitrary,margolin1997introduction, nobile1999stability, knupp2002reference, wells2005generation}.

The main goal of the location-based mesh movement approach is to directly
control the location of mesh points in particular regions of the computational
domain. A typical location-based method is the variational approach, which
relocates the mesh points by movements that are based on minimizing a functional
formulated to measure the difficulty or the error in the numerical
solution~\cite{huang2010adaptive}. Other location-based algorithms are based on
elliptic PDE descriptions which can be used to generate boundary-fitted meshes~\cite{thompson1974automatic,winslow1966numerical},  sometimes known as Winslow's approach~\cite{winslow1981adaptive}. Winslow's idea can be  generalized using a functional~\cite{brackbill1982adaptive}, which provides a combination of the mesh adaptivity, smoothness, and orthogonality conditions.

A number of articles consider mesh adaptation functionals including mechanical
models~\cite{jacquotte1988mechanical,jacquotte1992generation,jacquotte1992structured},
vector fields~\cite{knupp1995mesh}, a weighted Jacobian
matrix~\cite{knupp2002reference,knupp1996jacobian}, a matrix-valued diffusion
coefficient~\cite{huang1998high,cao1999study},  and the  equidistribution and
isotropy (or alignment conditions) presented in \cite{huang2001variational}. The
moving mesh PDE (MMPDE) method that we use has been developed by several authors~\cite{ren1992moving,huang1994moving,huang1994moving,huang1998high,huang1998moving,cao2005error}. Therein, the mesh movement is determined by a gradient flow equation,  and the functional plays the vital role.


To our knowledge, an adaptive moving mesh method, our method of choice, has not
been implemented for PDE-based modelling of pitting corrosion. In our moving
mesh method, the FEM provides the spatial discretization and our computational
framework is built upon the package MMPDElab by
Huang~\cite{huang2019introduction}. MMPDElab is a general  adaptive moving mesh
finite element solver for time dependent PDEs based on integration of the MMPDE.
The solver uses an alternating mesh and physical solution approach, and we used
the solver to achieve an adaptive moving mesh which provides sufficient mesh elements in and around the pit. 
We focus on the development of an appropriate mesh density function which implicitly and automatically determines an appropriate distribution of nodes as the pit evolves.


Our paper provides: 1) a proof of concept implementation of the moving mesh approach for pitting corrosion, 2) an adaptive solver for both single and multiple crystal directions, 3) the ability to handle single and multiple pits,  and 4) the ability to provide (provably) nonsingular quality evolving meshes in an automatic way.  The test material for demonstrating these techniques is 316 stainless steel using the parameters provided in~\cite{degiorgi2013numerical,kota2013microstructure}.
This paper is organized as follows. We discuss the preliminaries of the pitting
corrosion mechanism, the crystal orientation, the PDE-model, and an overview of
the moving mesh methodology in Section~\ref{sec:model}. The finite element
method approach for the physical PDE, the choice of mesh density function and
moving mesh parameters,
initial mesh generation, the boundary movement strategy, and the overall
alternating solution approach are discussed in
Section~\ref{sec:moving_mesh}.
Section~\ref{sec:numerical_results} is devoted to our numerical results and
validation.

\section{Model Problem}\label{sec:model}
In this section we detail our prototype model problem and adaptive solution strategy including the description of the domain, model PDEs and boundary conditions, the necessary crystallography, and an overview of the moving mesh strategy.

\subsection{PDE Model Equation}\label{sec:pde_model}
The transport of species in the solution can be modeled by the conservation of mass \cite{xie2015modeling,sharland1989finite}, which leads to the mathematical form
\begin{equation}\label{eq:pc02}
  \underbrace{\frac{\partial c_i}{\partial t}}_\text{Storage} = \underbrace{- \nabla N_i}_\text{Flux in - Flux out} + \underbrace{R_i}_\text{Generation},
\end{equation}
where $t$ is time,  $c_i$ is the concentration of the $i$th species, $N_i$ is the flux of the $i$th species, and $R_i$ is the rate of species generation due to chemical reactions. 
The ionic flux depends on the gradient of ion concentration,   electro-migration and convection (or flow in a liquid medium).
For each individual species $i$, the transport of the species in the electrolyte is described by the Nernst-Planck equation as 
\begin{equation}\label{eq:pc01}
  N_i = -D_i \nabla c_i - z_i F u_i c_i \nabla \varphi + c_i \textbf{v},
\end{equation}
where $D_i$ is the diffusion coefficient of the $i$th species, $u_i$ is the
mobility of the species, $z_i$ is the charge of the species, $\varphi$ is the
electric potential, $F$ is Faraday's constant, and \textbf{v} is the solvent
velocity \cite{bard1980electrochemical}.
Equation
(\ref{eq:pc01}) gives the flux of the species as a combination of three
contributing fluxes. The term $-D_i \nabla c_i$ describes the diffusive flux,
the term $- z_i F u_i c_i \nabla \varphi$ gives
the electro-migration flux and the term  $c_i \textbf{v}$ is the convection flux. Generally, the diffusion coefficient varies with the position of the species but we assume the diffusion coefficient is constant in our model.


While equation (\ref{eq:pc01}) represents the complete coupled behavior observed in the electrolyte domain during pitting corrosion, certain physically motivated assumptions are made in arriving at our simplified model:
a) the absence of gradients in the species concentration  due to the rapid mixing of the electrolyte;
b) the incompressibility of the solvent; and
c) the zero net production of
the reactants.

These assumptions simplify \eqref{eq:pc02} to the well-known Laplace equation.
In this case, the electrolyte potential can be found by solving (on the
electrolyte domain, $\Omega$, shown in Figure \ref{fig:model_domain})
\begin{minipage}{0.50\textwidth}
  \begin{equation}\label{eq:laplace}
    \nabla^2 \varphi = 0 \text{ in } \Omega,
  \end{equation}
\noindent  with the following boundary conditions
  \begin{equation}\label{eq:laplace_bc}
    \begin{split}
      \varphi  & = 0  \text{ on } \Gamma_1, \\
      \nabla \varphi\cdot \mathbf{n}  & = 0 {\text{ on }} \Gamma_2,\Gamma_3,\Gamma_4,\\
      \nabla \varphi\cdot \mathbf{n} & =\frac{i{(\varphi)}}{\sigma_c} \text{ on }\Gamma_p,
    \end{split}
  \end{equation}
\end{minipage}
~\hspace*{\fill}
\begin{minipage}{0.45\textwidth}
  \centering
\includegraphics[height=4.5cm]{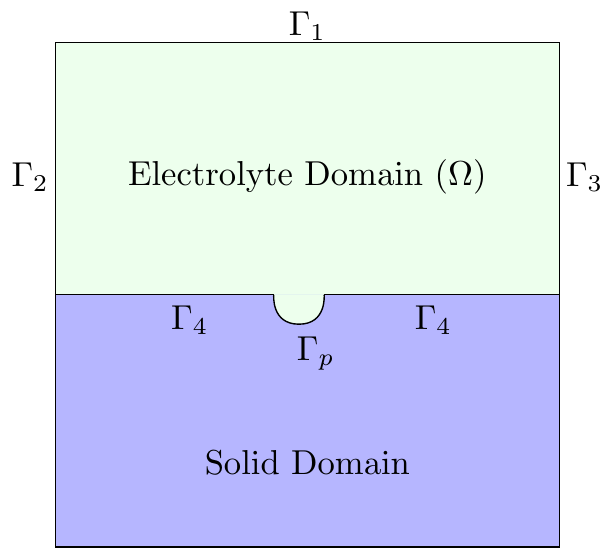}
\captionof{figure}{\label{fig:model_domain}The 2D computational domain.}
\end{minipage}\\
where $\nabla \varphi\cdot \mathbf{n} = \frac{\partial \varphi}{\partial
  \mathbf{n} }$, $\mathbf{n}$ is the (outward) unit normal vector, $i(\varphi)$ is the current density, $\sigma_{c}$ is the electrical conductivity of the electrolyte, $\Gamma_p$ is the pit boundary, and $\Gamma_1$, $\Gamma_2$, $\Gamma_3$  and $\Gamma_4$ are the top, left, right and bottom of the domain (excluding the pit boundaries), respectively. The boundary condition  $\frac{\partial \varphi}{\partial {\bf n}}  = 0 $ enforced on  $\Gamma_2$, $\Gamma_3$ and $\Gamma_4$ ensures there is no flow of ions across these boundaries. 
We denote the horizontal and vertical co-ordinates of the electrolyte region in
Figure \ref{fig:model_domain} by $x$ and $y$, respectively.

The current density is modelled by the Butler-Volmer relation
\begin{equation}\label{eq:bulter}
  i(\varphi) = z F A_s\cdot e^{\Big( \frac{z F ( V_{corr} + \alpha\eta_a)}{RT}  \Big)}, 
\end{equation}
where $\alpha$ is the transfer coefficient, $A_S$ is the material dissolution
affinity, $T$ is the temperature and $R$ is the universal gas constant \cite{bard2022electrochemical}.

The Butler-Volmer relation is used to describe the experimental data as a function of the applied over-potential
\[
  \eta_a = V_{app} - V_{corr} - \varphi,
\]
where $V_{app}$ and  $V_{corr}$ are the applied and the corrosion potentials, respectively~\cite{scheiner2009finite}.

As the metal corrodes, the pit boundary moves as the pit becomes larger.  In our
model the {new position of corrosion front}, $X_{\text{new}}$,  is computed from
the old position, $X_{\text{old}}$, by a simple time stepping procedure 
$$X_{\text{new}} = X_{\text{old}} + \Delta t V_{n}\mathbf{n},$$
where
$V_n$ is the  magnitude of normal velocity.
The magnitude of normal velocity at the corrosion interface (or the movement of the corrosion front) is  described using Faraday's law
\begin{equation}\label{eq:vn}
  V_{n} = \frac{i(\varphi)}{z F c_{solid}},
\end{equation}
where $c_{solid}$ is the atomic mass concentration of the metal and $z$ is the average charge number for the metal.  A list of parameters is given in Table \ref{tbl:list}.

\begin{table}
  \centering
  \begin{tabular}{@{}|clll|@{}} \hline
    Parameters          & : & Description                            & Value\\ \hline
    $z $                & : & Average charge number for the metal    & 2.19             \\
    $F $                & : & Faraday's constant                     & 96485 C/mol      \\
    $R $                & : & Universal gas constant                 & 8.315 J/(mol K)  \\
    $T $                & : & Temperature                            & 298.15 K         \\
    $V_{\text{corr}} $  & : & Mean corrosion potential (homogeneous) & -0.24 V          \\
    $V_{\text{app}} $   & : & Applied potential                      & -0.14 V          \\
    $A_{\text{diss}} $  & : & Dissolution affinity                   & 4 mol/cm$^2$s    \\
    $C_{\text{solid}} $ & : & Solid concentration                    & 143 mol/l        \\
    $\Delta t$          & : & Time step size                         &                  \\
    \hline
  \end{tabular}
  \caption{List of parameters used in the corrosion model.}\label{tbl:list}
\end{table}

\subsection{Crystal Orientation and Corrosion Potential}
The corrosion potential is the mathematical link between the etching effects of the electrolyte and the material undergoing pitting. This potential term is present in the Butler-Volmer equation \eqref{eq:bulter} and 
is an important parameter governing the velocity of the sides of the pit during
corrosion.  If the electrolyte etches the material homogeneously in all
directions then the crystal structure is not a variable in the
modeling and only one number is required for $V_{corr}$. However, if the
material is crystalline then the corrosion potential may vary dependent on the
particular crystallographic surface exposed to the electrolyte. Hence, a
connection between the Cartesian $(x,y,0)$ geometry used for defining the
computational domains as presented in Figure~\ref{fig:model_domain} and the
directions in the crystal, is needed.  In general, these two geometries will not
align since crystals can be rotated to lie along an infinite number of
directions and a transformation will be required to relate the two coordinate
systems.
Letting $\mathbf{n}_{CD}$ represent the outward unit normal in the crystal
coordinate system, and using the notation given in \cite{boisen2018mathematical,robertson2008imaging}, its relationship to
$\mathbf{n}$ is
\begin{equation}\label{eq:ncd_nxy}
  \mathbf{n}_{CD}=\mathbf{A}^{-1} \mathbf{n}.
\end{equation}
The matrix $\mathbf{A}^{-1}$ is defined by
\begin{equation}\label{eq:rotation}
  \mathbf{A}^{-1} = \left[
    \begin{array} {c|c|c}
      \          & \          & \          \\
      \mathbf{i} & \mathbf{j} & \mathbf{k} \\
      \          & \          & \          \\
    \end{array} \right],
\end{equation}
where $\mathbf{i}$, $\mathbf{j}$ and $\mathbf{k}$ are orthogonal unit column vectors in Cartesian space.  $\mathbf{k}$ is chosen as the desired zone axis of the crystal and $\mathbf{i}$ is chosen as the direction perpendicular to $\mathbf{k}$ which will be oriented along the $x$ axis in the computational domain.  Thus, the crystal can be rotated in any direction about the zone axis offering maximum flexibility in the problems that can be studied.  The third unit vector, $\mathbf{j}$, is found using the perpendicularity property of vector cross products
\begin{equation*}
  \mathbf{j}=\mathbf{k}\times \mathbf{i}.
\end{equation*}

In order to better see how to form the $\mathbf{A}^{-1}$ matrix, two examples
will be presented.  First, select a crystal orientation where the zone axis is
aligned along $\mathbf{k}=[001]$.  Next, choose the $\mathbf{i}=[100]$ direction
to be along the $x$ computational domain direction.  Performing the
$\mathbf{j}=\mathbf{k}\times \mathbf{i}$ cross product, it is found that
$\mathbf{j}=[010]$.  Hence, the transformation matrix is
\begin{equation*}
  \mathbf{A}^{-1}_{001} = \left[
    \begin{array} {ccc}
      1 & 0 & 0 \\
      0 & 1 & 0 \\
      0 & 0 & 1 \\
    \end{array} \right].
\end{equation*}
In this case, the transformation matrix is the identity matrix and the crystal coordinate system is the same as the computational domain coordinate system, ie.\ $\mathbf{n}_{CD}=\mathbf{n}$.  For the second example, choose the zone axis to be along the [101] direction and the $[\bar{1} 0 1]$ crystal direction to be along the $x$ computational domain direction.  This leads to $\mathbf{j}=[010]$.  After normalizing the vectors, the transformation matrix in this case is
\begin{equation*}
  \mathbf{A}^{-1}_{110} = \left[
    \begin{array} {ccc}
      \frac{-1}{\sqrt{2}} & 0 & \frac{1}{\sqrt{2}} \\
      0                   & 1 & 0                  \\
      \frac{1}{\sqrt{2}}  & 0 & \frac{1}{\sqrt{2}} \\
    \end{array} \right].
\end{equation*}

\begin{figure}[ht!]
  \begin{center}
    \includegraphics[width=4.5in]{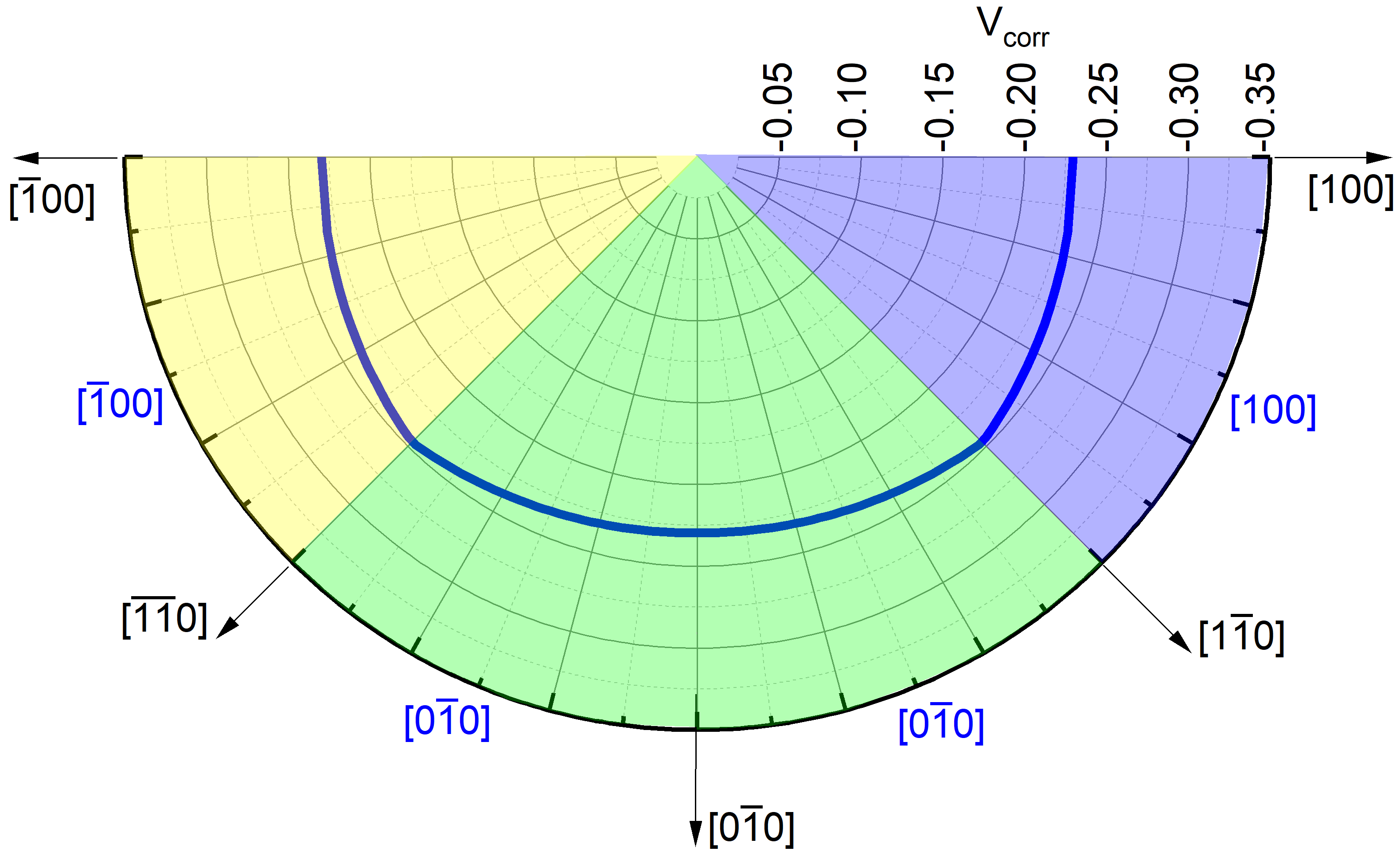}
    \caption{A plot of $V_{corr}$ as a function of location around the pit boundary.  The blue, green and yellow sections of the pit edge require different $<\!001\!>$ vectors for use in equation \eqref{eq:vcoor}.}
    \label{fig:fig2}
  \end{center}
\end{figure}

The next step is to define a corrosion potential for each crystallographic
direction.  Unfortunately, experimental data giving the corrosion potential as a
function of crystallographic surface is usually not available, hence we will adopt a similar form of semi-empirical potential for 316 stainless steel 
\begin{equation}\label{eq:vcoor}
  V_{corr} = k-s\left[1-\left(\left\langle 001\right\rangle \cdot \mathbf{n}_{CD}\right)_{\textrm{max}}\right],
\end{equation}
where $k=-0.2297$ and $s=0.054$ gives a 10\% difference between the maximum and
minimum $V_{corr}$ values, that is between the [001] and [111] crystal planes as
used by DeGiorgi et al.\ in \cite{kota2013microstructure}.  
We write  $\left(\left\langle 001\right\rangle \cdot
  \mathbf{n}_{CD}\right)_{\textrm{max}}$, rather than  $\left[ 001\right] \cdot
\mathbf{n}_{CD}$ as used in \cite{kota2013microstructure}.  In our case, $\left\langle 001\right\rangle$ represents any one of the six cryptographically equivalent $[001]$, $[\bar{1}00]$, $[010]$, $[0\bar{1}0]$, $[001]$, $[00\bar{1}]$ directions that maximizes the dot product with the crystal direction normal vector.  Maximizing this dot product minimizes the angle between the normal vector and the particular $\left\langle 001\right\rangle$ vector so that an equivalent result to the standard stereographic triangle is obtained.  An example of this procedure is provided in Figure \ref{fig:fig2} where a single crystal has been oriented along the [001] zone axis.  The semicircle represents the pit boundary, the black lines and arrows the outward pit normal vectors, and the heavy blue line is the value of $V_{corr}$ as a function of position along the edge of the pit.  Highlighted by the three colours are sections along the pit boundary that have a different $\left\langle001\right\rangle$ vector for use in equation (\ref{eq:vcoor}).  These vectors are presented in blue text and are $[001]$, $[0\bar{1}0]$ and $[\bar{1}00]$ for the right, centre and left sections of the pit, respectively.




\subsection{Overview of the Moving Mesh Strategy}
The basic idea of the moving mesh method is to automatically redistribute a
fixed number of nodes where additional accuracy is required. The mesh moves or
evolves automatically as the solution or domain evolves and is obtained by
solving a MMPDE. This MMPDE depends on a mesh density or
monitor function, which is large where the mesh density is needed to be large.
The mesh density function is often chosen to depend on variations or errors in
the solution of the physical PDE or is chosen by geometrical considerations (as in this paper).

In 1D, the equidistribution principle (EP) is used to derive the moving mesh
system. The EP uses a mesh density function $\boldsymbol{\rho} =
\boldsymbol{\rho}(x) >0$ which is to be distributed evenly among the mesh elements in the domain. Given an integer $N > 1$,  the continuous and bounded function $\boldsymbol{\rho}$ on $[a,b]$ is evenly distributed on the mesh $\mathcal{T}_h = a = x_0<x_1<\cdots <x_N =b$, if
\begin{equation}\label{eq:ep}
  \int_{x_0}^{x_1} \boldsymbol{\rho} (x)dx =\int_{x_1}^{x_2} \boldsymbol{\rho} (x)dx  = \cdots
  =\int_{x_{N-1}}^{x_{N}} \boldsymbol{\rho} (x)dx.
\end{equation}
A mesh $\mathcal{T}_{h}$ is called an equidistributing mesh if the mesh satisfies the equidistribution principle.

The physical problem is assumed to require a non-uniform $x$-coordinate,
$x\in\Omega$.  This physical coordinate, $x$,  is a mapping of the
computational $\xi$-coordinate where $\xi \in \Omega_c = [0, 1]$, and $x(0) = a$ and $x(1) = b$, if $\Omega=[a,b]$. We attempt to generate a physical mesh $\mathcal{T}_h$ using a mesh transformation $x = x(\xi): \Omega_c \rightarrow \Omega$ and a uniform mesh in the $\xi$-coordinate
\[
  \xi_i = \frac{i}{N}, \quad i= 0, 1,\ldots,N.
\]
The equidistribution principle \eqref{eq:ep} can then be written as
\[\begin{split}
    \int_{a}^{x_{i}} \boldsymbol{\rho} (\tilde{x})d\tilde{x} &= \frac{i}{N}\int_{a}^{b} \boldsymbol{\rho} (\tilde{x})d\tilde{x}\\
    &=  \frac{i}{N} \boldsymbol{\sigma}, \quad i= 0, 1,\ldots,N,
  \end{split}
\]
where $ \boldsymbol{\sigma} = \int_{a}^{b} \boldsymbol{\rho}
(\tilde{x})d\tilde{x}$. The function $\int_{a}^{x} \boldsymbol{\rho}
(\tilde{x})d\tilde{x}$ is strictly monotonically increasing if
$\boldsymbol{\rho}>0$, therefore each $x_i$ is unique. Here
$\boldsymbol{\sigma}$ and $\frac{i}{N} \boldsymbol{\sigma}$ are the total error
and average error in the approximating solution, respectively.

Using the mesh transformation we have
\[
  \int_{a}^{x(\xi_{i})} \boldsymbol{\rho} (\tilde{x})d\tilde{x} = \xi_i \boldsymbol{\sigma}, \quad i= 0, 1,\ldots,N,
\]
and the continuous version is given by 
\begin{equation}\label{eq:ep_cont}
  \int_{a}^{x(\xi)} \boldsymbol{\rho} (\tilde{x})d\tilde{x} = \xi \boldsymbol{\sigma}, \quad \forall \xi \in \Omega_c.
\end{equation}
The continuous mapping $x=x(\xi)$ is called a \textit{equidistributing
  coordinate transformation} for $\boldsymbol{\rho}$ if it satisfies
relation~\eqref{eq:ep_cont}. Differentiating (\ref{eq:ep_cont}) with respect to $\xi$ gives
\begin{equation}\label{eq:equid_1}
  \boldsymbol{\rho} (x)\frac{d x}{d\xi} = \boldsymbol{\sigma}.
\end{equation}
Equation~\eqref{eq:equid_1} indicates that $\frac{d x}{d\xi} $ is small when  $\boldsymbol{\rho}$ is large. Again, differentiating with respect to $\xi$ gives
\begin{equation}\label{eq:equid_2}
  \frac{d}{d\xi}\Big(\boldsymbol{\rho} (x)\frac{d x}{d\xi}\Big) = 0,
\end{equation}
\noindent with the boundary conditions
\begin{equation}\label{eq:equid_bc1}
  x(0) = a, \quad x(1) = b.
\end{equation}
This is a nonlinear boundary value problem (BVP) for the required mesh
transformation and physical mesh. The mesh and physical solution on that mesh is
determined by solving this BVP and the physical PDE as a coupled system.

In higher dimensions, in order to describe the equidistribution and alignment
conditions at the discrete level, we consider a  mesh $ \mathcal{T}_h$ of $N$
triangular elements with $N_v$ vertices in the physical domain $\Omega \in
\mathbf{R}^d$ $(d \ge 1)$.  Furthermore, we consider an invertible affine mapping $F_K : \hat{K} \rightarrow K$ and its Jacobian matrix, $F_K^{'}$, where $\hat{K}$ is the reference or master element for a physical element $K$ in $\mathcal{T}_h$. 
Assume that a metric tensor (or a monitor function) $\mathbb{M} = \mathbb{M}(x)$ is given on $\Omega$ which determines the shape, size and orientation of mesh elements of the domain $\Omega$. Generally, a mesh is uniform if all of its elements have the same size and is similar to a reference element $\hat{K}$. So, the main idea of the MMPDE method is to view any adaptive mesh  $\mathcal{T}_h$ as a uniform mesh in the metric $\mathbb{M}$.

The requirements of the equidistribution and alignment in higher dimensions can be expressed mathematically at the discrete level~\cite{HR2011} as
\[
  |K| \sqrt{\texttt{det}(\mathbb{M}_K)} = \frac{\sigma_h}{N}, \forall K\in \mathcal{T}_h,
\]
and
\[
  \frac{1}{n}tr((F_K^{'})^T\mathbb{M}_K F_K^{'}) = \texttt{det}((F_K^{'})^T\mathbb{M}_K F_K^{'})^{\frac{1}{d}}, \forall k \in \mathcal{T}_h,
\]
where $|K|$ is the volume of $K$ and $\sigma = \sum_{K \in \mathcal{T}_h}|K| \sqrt{\texttt{det}(\mathbb{M}_k)}$. 


A standard choice bases the metric tensor $\mathbb{M}$
on the approximate Hessian of the solution.
This choice of $\mathbb{M}$ is known to be optimal with respect to the $L_2$
norm of the linear interpolation error \cite{huang2005metric}.
Here we focus on the evolution of the pit geometry and choose $\mathbb{M}$
through geometrical considerations (see Section \ref{sec:meshdensity}).

A discrete functional associated with the equidistribution and alignment
conditions is given by
\begin{equation}\label{eq_energy_fun}
  I[\mathcal{T}_h] =  \sum_{K \in \mathcal{T}_h}|K|\  \sqrt[]{\text{det}(\mathbb{M}_K)}\ \Big[ \theta\ \Big(\text{tr}(\mathbb{J}\mathbb{M}_{K}^{-1}\mathbb{J}^{T})\Big)^{\frac{d\gamma}{2}} +
    (1-2\theta)
    {d^\frac{d\gamma}{2}}{}\Big(\frac{\text{det}(\mathbb{J})}{\sqrt{\text{det}(\mathbb{M}_K)}}\Big)^\gamma\Big],
\end{equation}
where $\mathbb{J} =   (F_K^{'})^{-1}.$ Minimizing this functional $I[\mathcal{T}_h]$ approximately satisfies the equidistribution and alignment conditions \citep{huang2001variational}. The value of the parameters $\theta = \frac{1}{3}$, and  $\gamma = \frac{3}{2}$ are used for our numerical experiments.

The MMPDE moving mesh equation can then be defined as the (modified)
gradient system (or gradient flow equation) for the energy functional, i.e.,
\begin{equation}\label{eq:gradientflow}
  \frac{d {\bf x}_i}{dt} = - \frac{P_i}{\tau} \frac{\partial I[\mathcal{T}_h]}{\partial {\bf x}_i}, \quad i = 1,2,\ldots, N_v,  \quad t\in (t_n, t_{n+1}],
\end{equation}
where $P_i = \text{det}(\mathbb{M}_i)^{\frac{1}{d+2}}$ is a scalar function used
to ensure
the mesh equation has invariance properties and $\tau$ is a positive parameter used to adjust the response time of
mesh movement to the change in $\mathbb{M}$. A smaller value of $\tau$ provides
 a faster response. 

\section{The Numerical Implementation}\label{sec:moving_mesh}
This section describes the details of the adaptive MMPDE strategy used to solve the PDE pitting corrosion model using a customized version of MMPDELab~\cite{huang2019introduction}.
\subsection{Discretization Of The Physical PDE}
MMPDElab requires the user to specify the physical PDE in weak form, where 
the strong form of our model problem is given in equations \eqref{eq:laplace} and \eqref{eq:laplace_bc}.
Let $V$ be the trial space, chosen in this case as
\[
  V = \{ v \in H^{1}(\Omega(t)) : v = 0 \text{ on }          \Gamma_1\} \subset  H_{}^{1} (\Omega(t)),
\]
where $H_{}^{1}(\Omega(t))$ is, roughly speaking, the function space whose
members, and their first derivatives,  are square integrable (see, for example,
\cite{Adams_David_R_and_Hedberg_Lars_Inge1999-oy} for details).
At any time $t$ the weak form is constructed as follows: find $\varphi \in V$ such that
\begin{equation}\label{eq:variational}
  \int_{\Omega(t)}\nabla \varphi\cdot \nabla v d\Omega = \int_{ \Gamma_p(t)}v \frac{i{(\varphi)}}{\sigma_c} ds,  \quad \forall v \in V,
\end{equation}
where $\Gamma_p(t)$ is the boundary of the pit at time $t$.
Here $V_h$ denotes a finite dimensional subspace of $V$   spanned by a
collection of finitely many basis functions (often associated with a mesh).  We
discretize the weak  form~\eqref{eq:variational} to find a solution in the
discrete trial space. The discrete FEM solution is then found by finding $\varphi_h \in V_h \subset V$, such that
\begin{equation}\label{eq:fem}
  \int_{\Omega(t)}\nabla \varphi_h \cdot \nabla v_h d\Omega = \int_{ \Gamma_p(t)} \frac{i{(\varphi_h)}}{\sigma_c} v_h ds,  \quad \forall v_h \in V.
\end{equation}

We can solve the discrete variational problem \eqref{eq:fem} in the following way. First, introduce $\big\{\phi_j\big\}_{j=1}^{N}$ as a basis for $V_h$ and $V$. Let $\varphi_h \in V_h$ be a linear combination of the basis functions $\phi_j, j=1,2,\ldots,N$, with coefficients $\tilde{\varphi}_j$,  given by
\begin{equation}\label{eq:basis}
  \varphi_h = \sum_{j=1}^{N} \tilde{\varphi}_j \phi_j.
\end{equation}
Considering $v = {\phi}_k$, for $k=1,2,\ldots,N$, and using relation~\eqref{eq:basis} gives
\begin{equation*}\label{eq:fem3}
  \sum_{j=1}^{N} \tilde{\varphi}_j\int_{\Omega}\nabla \phi_j. \nabla  {\phi}_k d\Omega -  \frac{1}{\sigma_c} \int_{ \Gamma_p} i{\Big(\sum_{j=1}^{N} \tilde{\varphi}_j \phi_j\Big)} {\phi}_k ds  =0, \quad k=1,2,\ldots,N. \\
\end{equation*}
For each time, $t$, this is a system of non-linear equations which is solved using Newton's method.

%

\subsection{The Choice Of The Mesh Density Function}\label{sec:meshdensity}
The appropriate specification of the mesh density tensor is crucial --- it controls how the mesh automatically adapts to the changing solution features.

As mentioned, for problems with fixed domain boundaries,  the Hessian based monitor function is an often used, general purpose, driver of the adaptive mesh.   Here, however, we wish to ensure sufficient resolution of the pit geometry.
To ensure a sufficient number of elements in the evolving pit and near the pit
boundary, we use a modified version of MacKenzie's  distance-based monitor function
\begin{equation} \label{eq_macKenzie_mf_2}
  \mathbb{M}_{K}(x,y) = \left(1+\frac{\mu_1}{\sqrt{\mu_2^2 d^2 + 1}}\right) I,
\end{equation}
where
$$d(x,y) = \min_{p}{\|(x, y)-(x_p, y_p)\|},$$
and ($x_p,y_p$) denotes any point on the boundary of the pit, $\Gamma_{pit}$,
cf.\ \cite{beckett2001moving}.
At any point $(x,y)\in\Omega$ the value of the monitor function involves the
minimum distance, measured in the two-norm $\|\cdot\|$, from $(x,y)$ to any point on the pit boundary.   The reciprocal of $\mathbb{M}_{K}$ indicates that $\mathbb{M}_K$ will be largest in $(x,y)$ regions where the distance to the pit boundary is the smallest, and hence the mesh spacing will be automatically smaller in these regions.
The parameter $\mu_1$ controls the minimum
mesh spacing whereas $\mu_2$ (and $\tau$) will control the rate at which mesh clustering occurs during the integration of the MMPDE \cite{mackenzie2000numerical}.

To understand the effect of the $\mu_1$ and $\mu_2$ parameters we consider a simple experiment, simulating the evolution of a single pit in a homogeneous material.   A quality initial mesh, generated using the process outlined in Section
\ref{sec:initialmeshgeneration}
is used for this experiment.


\begin{figure}[ht!]
\centering
  \begin{subfigure}{.40\textwidth}
    \includegraphics[width=.98\linewidth,height=.65\linewidth]{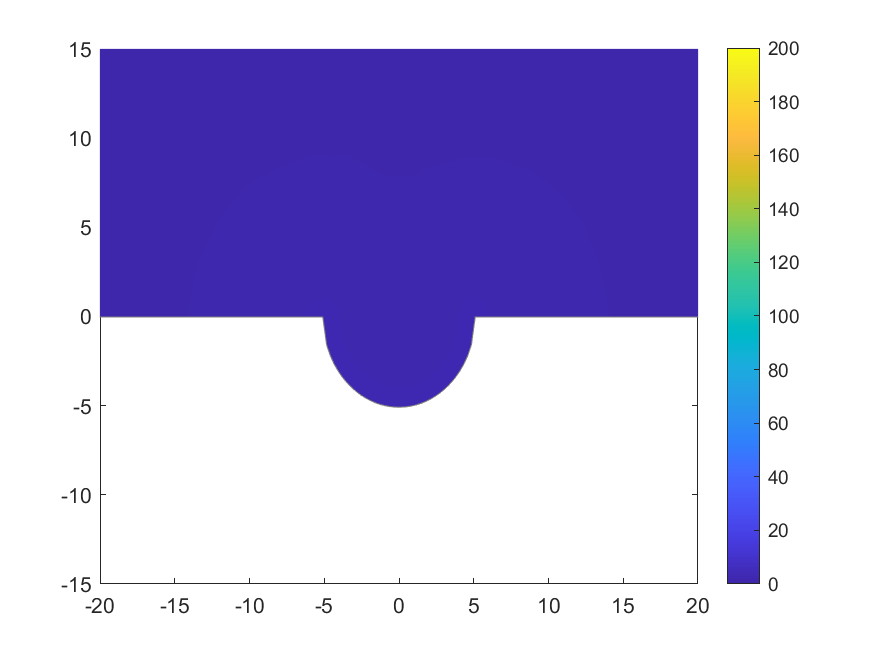}
    \caption{Monitor function values with $\mu_1 = 1$.}
  \end{subfigure}
  \begin{subfigure}{.40\textwidth}
    \includegraphics[width=.98\linewidth,height=.65\linewidth]{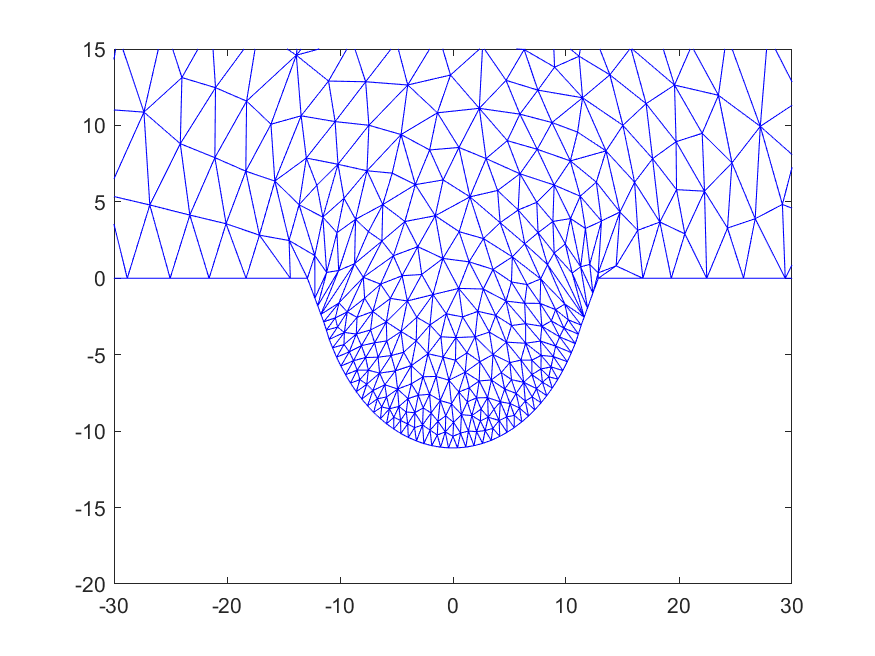}
    \caption{The mesh at $t=120$ s with $\mu_1 = 1$.}
  \end{subfigure}

  \begin{subfigure}{.40\textwidth}
    \includegraphics[width=.98\linewidth,height=.65\linewidth]{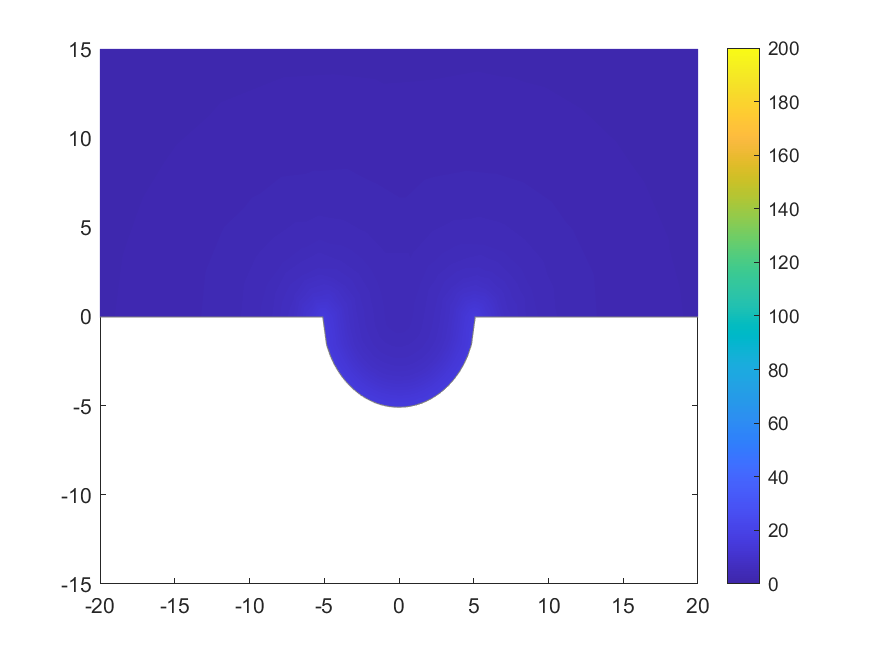}
    \caption{Monitor function values with $\mu_1 = 10$.}
  \end{subfigure}
  \begin{subfigure}{.40\textwidth}
    \includegraphics[width=.98\linewidth,height=.65\linewidth]{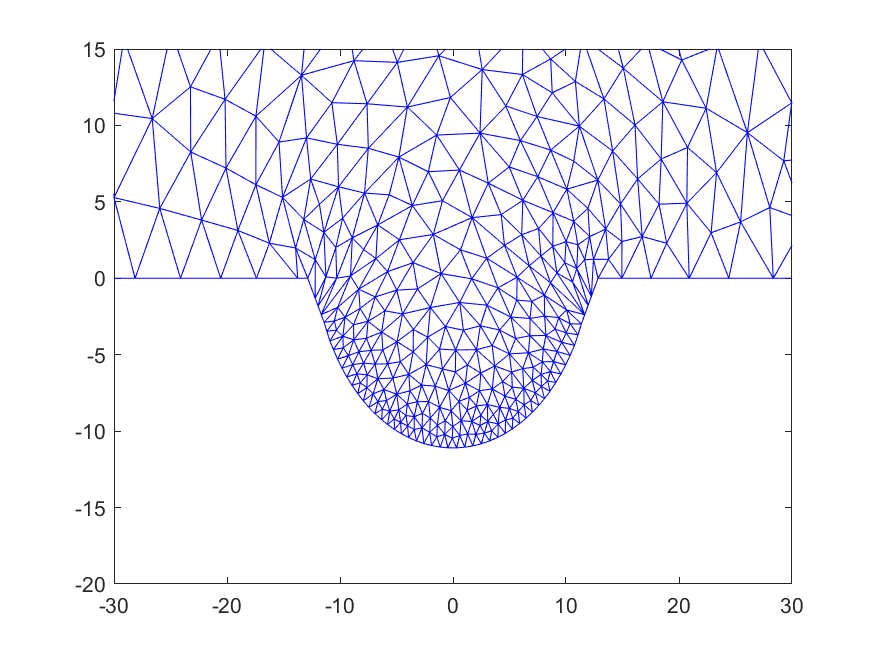}
    \caption{The mesh at $t=120$ s with $\mu_1 = 10$.}
  \end{subfigure}

  \begin{subfigure}{.40\textwidth}
    \includegraphics[width=.98\linewidth,height=.65\linewidth]{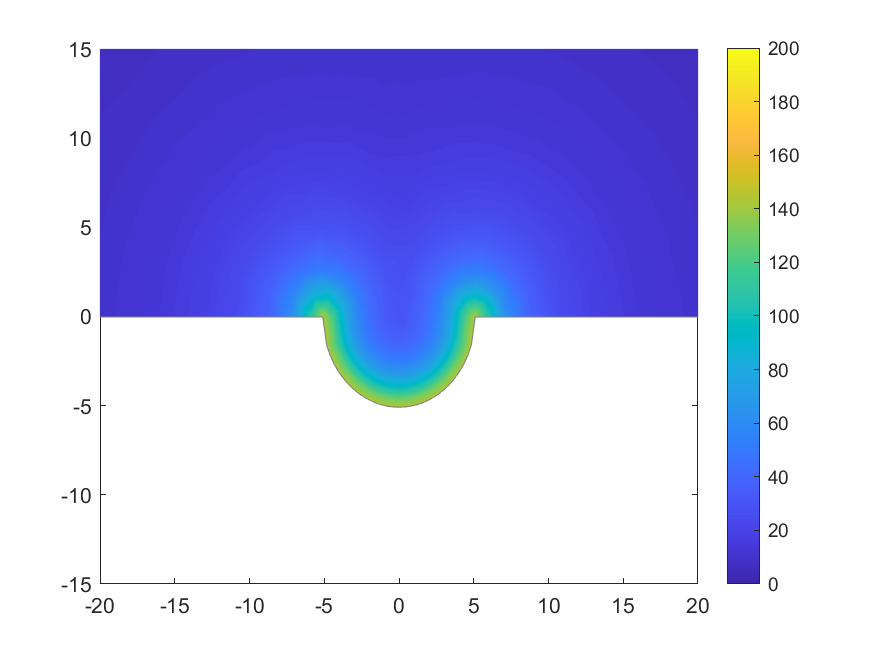}
    \caption{Monitor function values with $\mu_1 = 100$.}
  \end{subfigure}
  \begin{subfigure}{.40\textwidth}
    \includegraphics[width=.98\linewidth,height=.65\linewidth]{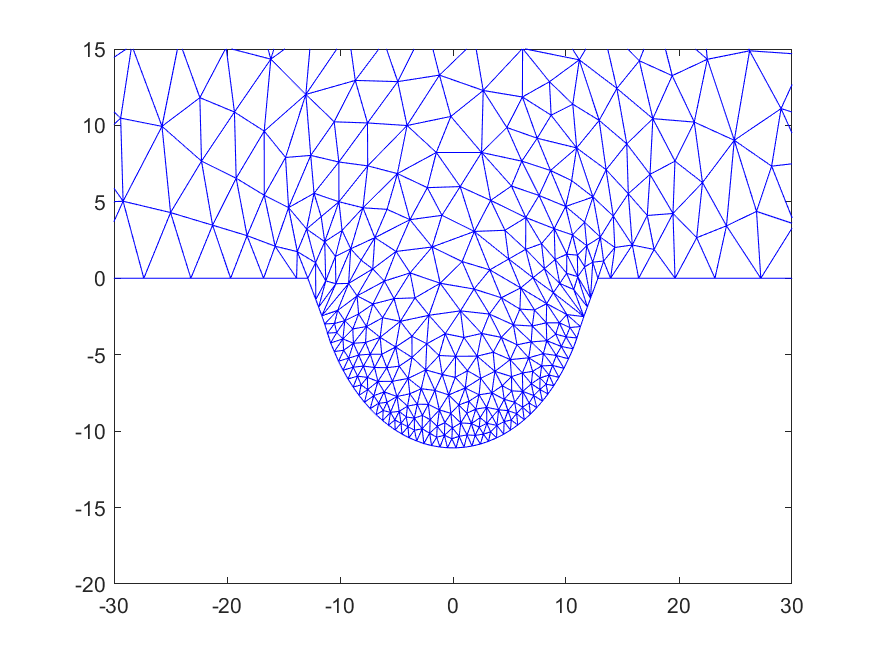}
    \caption{The mesh at $t=120$ s with $\mu_1 = 100$.}
  \end{subfigure}
  \caption{Effect of $\mu_1$ on the mesh at $t=120s$ for the simulation of a pit in a homogeneous material with the monitor function \eqref{eq_macKenzie_mf_2} and $\mu_2=1$.}\label{fig:effect-of-mu1}
\end{figure}

We begin by fixing $\mu_2$, at $\mu_2=1$, and consider the effect of increasing
$\mu_1$.  The plots on the left of each row in Figure \ref{fig:effect-of-mu1}
show a representative mesh density function (computed at $t=2s$), while the
computed mesh after $t=120s$ is shown on the right.  The results show that
increasing $\mu_1$ leads to a monitor function which is (relatively) larger near
the pit boundary and hence gives smaller grid spacings near the pit boundary.
The value $\mu_1=100$ provides a balance between increased mesh density near the pit boundary and sufficient resolution throughout the rest of the computational domain.

\begin{figure}[ht!]
\centering
  \begin{subfigure}{.40\textwidth}
    \includegraphics[width=.98\linewidth,height=.65\linewidth]{./figures/effect_of_mu/mu1_100_mu2_1/npit1_monitor_fun_contour_t2.png}
    \caption{Monitor function values using $\mu_2=1$.}
  \end{subfigure}
  \begin{subfigure}{.40\textwidth}
    \includegraphics[width=.98\linewidth,height=.65\linewidth]{./figures/effect_of_mu/mu1_100_mu2_1/npit1_mesh_t120.png}
    \caption{The mesh at $t=120$ s using $\mu_2=1$.}
  \end{subfigure}

  \begin{subfigure}{.40\textwidth}
    \includegraphics[width=.98\linewidth,height=.65\linewidth]{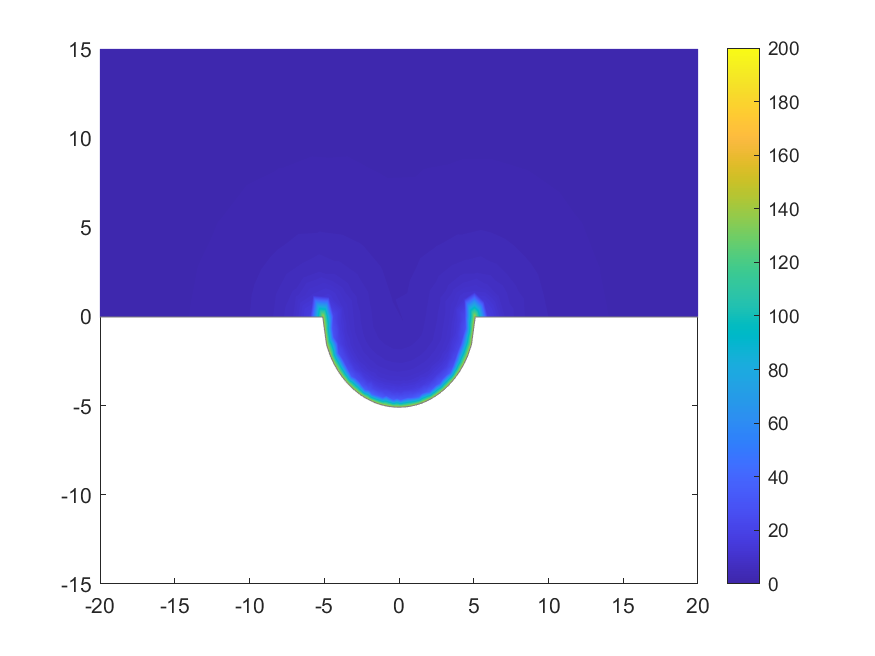}
    \caption{Monitor function values using $\mu_2=10$.}
  \end{subfigure}
  \begin{subfigure}{.40\textwidth}
    \includegraphics[width=.98\linewidth,height=.65\linewidth]{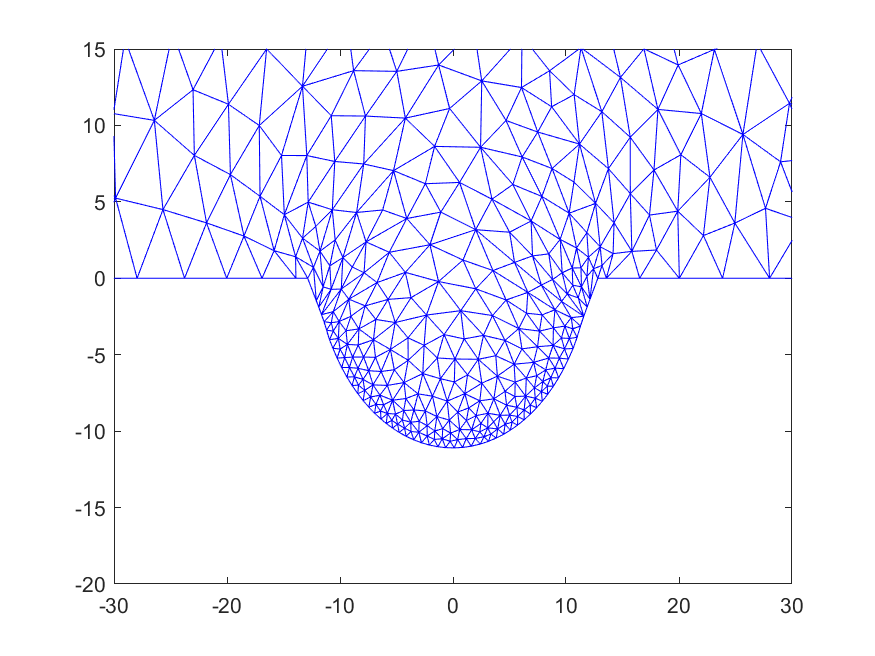}
    \caption{The mesh at $t=120$ s using $\mu_2=10$.}
  \end{subfigure}

  \begin{subfigure}{.40\textwidth}
    \includegraphics[width=.98\linewidth,height=.65\linewidth]{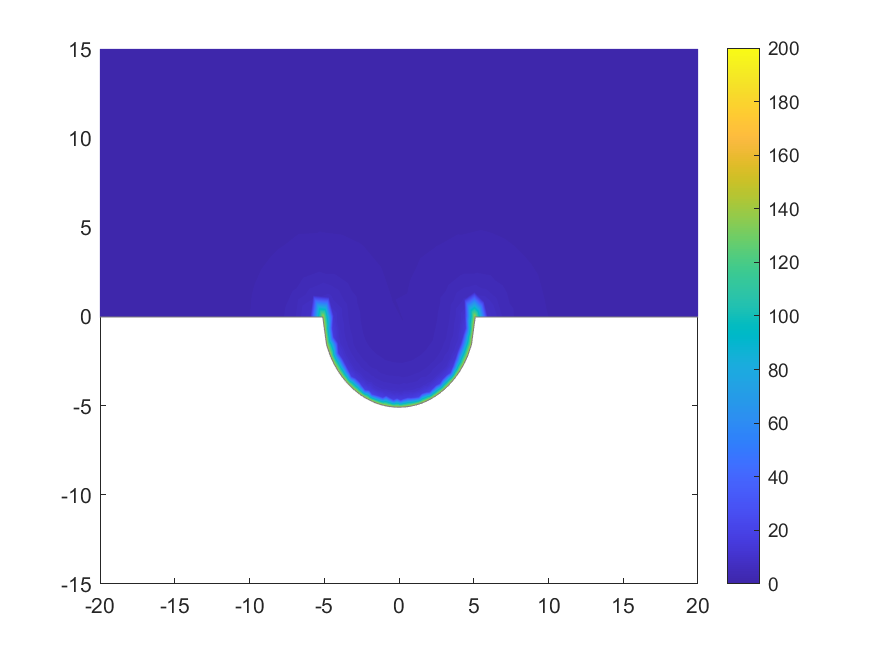}
    \caption{Monitor function value using $\mu_2=20$.}
  \end{subfigure}
  \begin{subfigure}{.40\textwidth}
    \includegraphics[width=.98\linewidth,height=.65\linewidth]{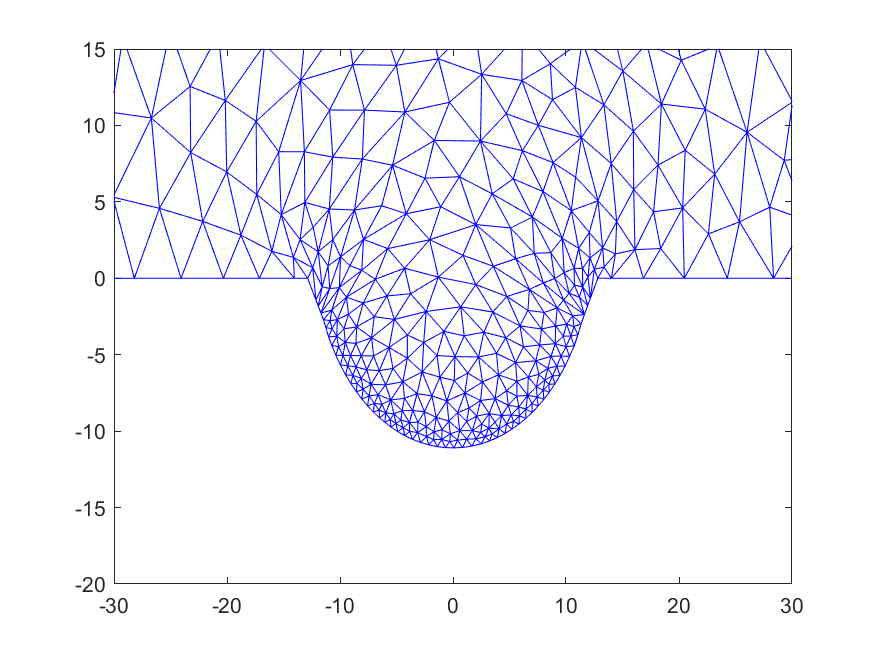}
    \caption{The mesh at $t=120$ s using $\mu_2=20$.}
  \end{subfigure}
   \caption{Effect of $\mu_2$ on the mesh at $t=120s$ during the simulation of a pit in a homogeneous material with the monitor function \eqref{eq_macKenzie_mf_2} and using $\mu_1=100$.}
  \label{fig:effect-of-mu2}
\end{figure}

Mackenzie \cite{mackenzie2000numerical} reports that increasing $\mu_2$ reduces
the spatial extent of node clustering near the pit boundary.   To explore this
we fix $\mu_1=100$ and vary $\mu_2$,  recording representative mesh density
function values and the final mesh obtained for the propogation of a homogeneous
pit is shown in Figure \ref{fig:effect-of-mu2}.
The figure shows that increasing $\mu_2$ to $20$ is       better able to keep the mesh focused on the feature of    interest (the pit boundary in this case),        consistent with the general findings in                   \cite{mackenzie2000numerical}.

\subsection{Initial Mesh Generation}\label{sec:initialmeshgeneration}

The numerical simulations in this paper require an initial pit geometry and an initial spatial grid. With a prescription of the spatial domain many software platforms provide tools known as mesh generators for this purpose, for example {\bf initmesh} in Matlab or the {\bf mesh node} in COMSOL.  These tools require a description of the domain boundary and then generate a mesh subject to constraints on the mesh size (or number of nodes) and aspect ratios of the mesh elements.

In most cases the solutions found on the initial meshes generated by these approaches will not be optimal --- for example, there is no guarantee that the error in the numerical solution will be minimized.
There are alternatives for initial mesh generation that involve the MMPDE approach considered in this paper, and hence are consistent with the technique used for all subsequent time steps.

One technique for initial mesh tuning begins by using the simple mesh generators mentioned above to
find a (nearly) uniform mesh.  The physical problem defined by equations
(\ref{eq:laplace})  and (\ref{eq:laplace_bc}) is then solved on this mesh to
give an initial potential.   Using this initial potential and its associated
mesh density function $\mathbb{M}_K$,  the gradient flow equation
(\ref{eq:gradientflow}) can be solved to a steady state (alternating its
solution with physical solves).  The result is a mesh which minimizes the
discrete functional (\ref{eq_energy_fun}), equidistributing the initial
potential over the initial computational domain.   Should a more sophisticated
non-uniform mesh generator be available, then an initial non--uniform mesh can
be smoothed in the same manner.  In practice, equation (\ref{eq:gradientflow})
may not solved to a steady state.  Instead, equation (\ref{eq:gradientflow}) can
be integrated for a specified number of time steps, or can be integrated until a
specified difference between two meshes is found.   We will call this mesh
smoothing. The number of steps required to reach a suitable approximation of the steady state is a function of the physical solution,  the number of mesh nodes, and the mesh density parameters.
This idea of using the MMPDE to tune the initial mesh has the added benefit of giving a mesh which has the same properties as all subsequent meshes,   while using the same code base as the rest of the simulation.
We note that this process can also be used to provide small scale mesh smoothing
during the solution of the  moving boundary problem.  This is particularly useful should the  pit boundary movement be large or if the pit boundary     movement induces a discontinuous change in the geometry   (during pit merging for example).

\begin{figure}[ht!]
  \begin{subfigure}{.5\textwidth}
    \centering
    \includegraphics[width=.98\linewidth, height=.65\linewidth]{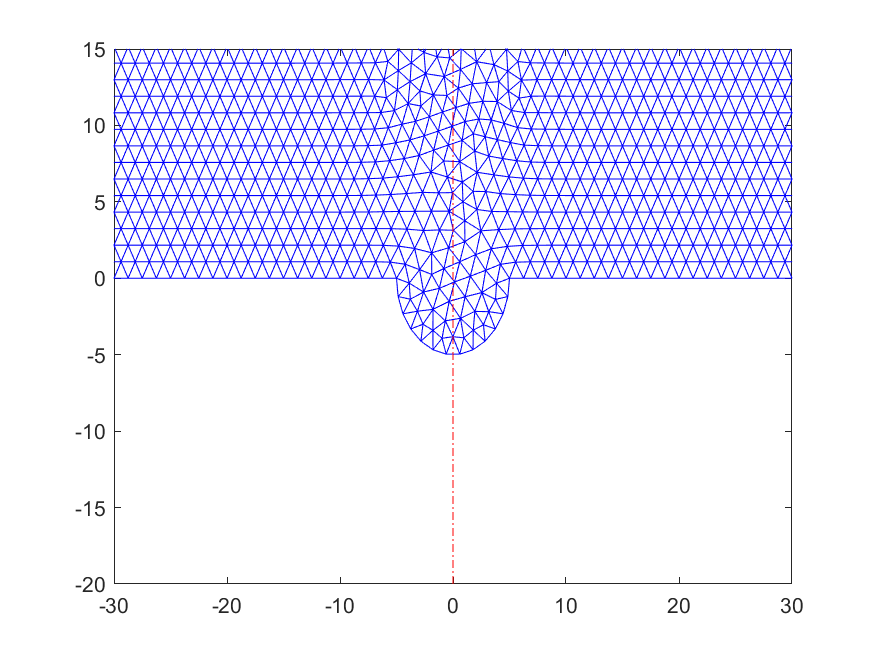} 
    \caption{An uniform initial mesh.}
    \label{fig:uniforminitial}
  \end{subfigure}
  \begin{subfigure}{.5\textwidth}
    \centering
    \includegraphics[width=.98\linewidth, height=.65\linewidth]{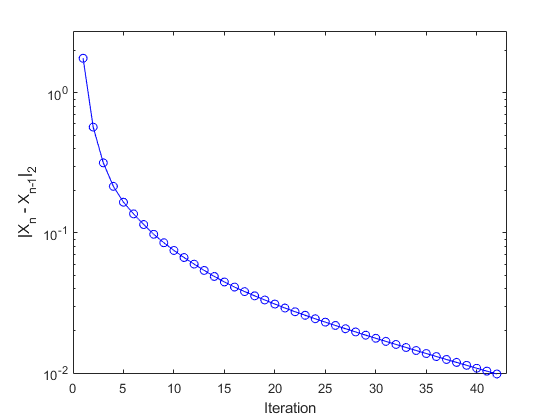} 
    \caption{The convergence of the initial mesh smoothing.}
    \label{fig:uniformconvergence}
  \end{subfigure}

  \begin{subfigure}{.5\textwidth}
    \centering
    \includegraphics[width=.98\linewidth, height=.65\linewidth]{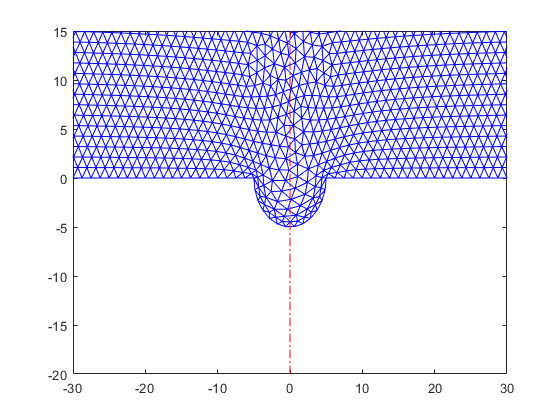} 
    \caption{The smoothed initial mesh.}
    \label{fig:smoothedinitial}
  \end{subfigure}
  \begin{subfigure}{.5\textwidth}
    \centering
    \includegraphics[width=.98\linewidth, height=.65\linewidth]{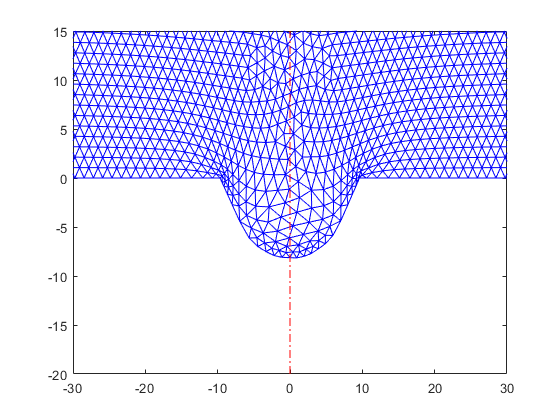} 
    \caption{The mesh after $60$ s.}
    \label{fig:finalmesh2}
  \end{subfigure}

  \caption{(a) Uniform initial mesh, (b) convergence of the mesh smoothing
    process (c), initial mesh after mesh smoothing, and (d) the mesh after 60 s using the monitor function \eqref{eq_macKenzie_mf_2}.}
  \label{fig:uniform-mesh-over-time}
\end{figure}

In Figure~\ref{fig:uniform-mesh-over-time} we show the evolution of the pit
geometry in a homogeneous material.   The moving mesh method is implemented with
the monitor function~\eqref{eq_macKenzie_mf_2} starting from the initial uniform
mesh shown in Figure \ref{fig:uniforminitial}.  The non-uniform mesh resulting from the
solution of the gradient flow equation, shown in Figure \ref{fig:smoothedinitial}, 
  has successfully concentrated the mesh elements in
the initial pit and near the initial pit boundary.  The convergence to a steady
state solution of the gradient flow equation is shown in Figure
\ref{fig:uniformconvergence}. This resulting non-uniform
mesh is now an appropriate initial mesh to use to  evolve the pitting corrosion
problem forward in time.    It is important to stress that no hand-tuning of the
initial mesh is necessary, it is generated automatically during the smoothing
process based on the characteristics of the chosen mesh density function.  The
final mesh after $t=60$ s is given in Figure \ref{fig:finalmesh2}.

\begin{figure}[ht!]
  \begin{subfigure}{.5\textwidth}
    \centering
    \includegraphics[width=.98\linewidth]{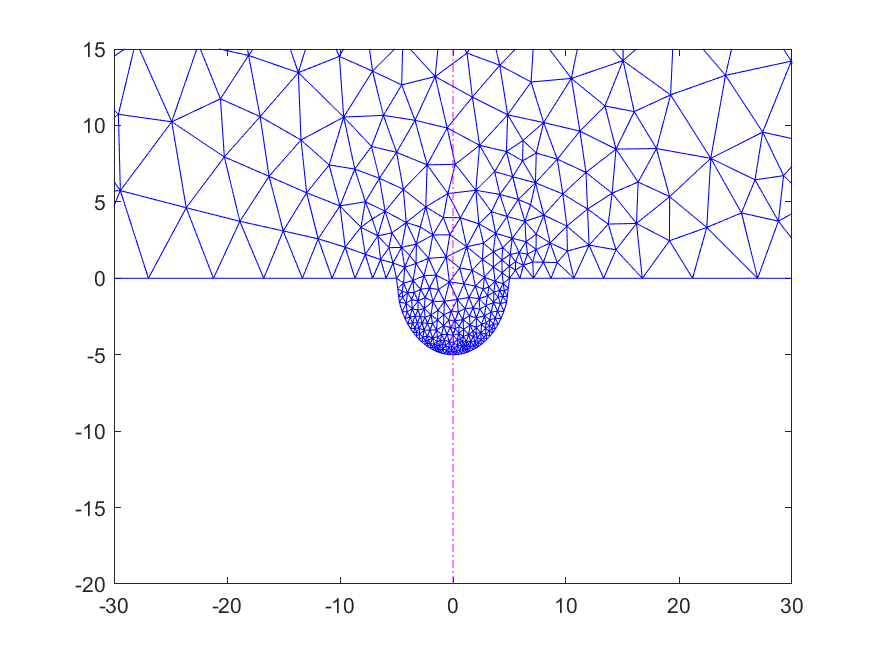}
    \caption{A nonuniform initial mesh.}\label{fig:nonuniform}
  \end{subfigure}
  \begin{subfigure}{.5\textwidth}
    \centering
    \includegraphics[width=.98\linewidth]{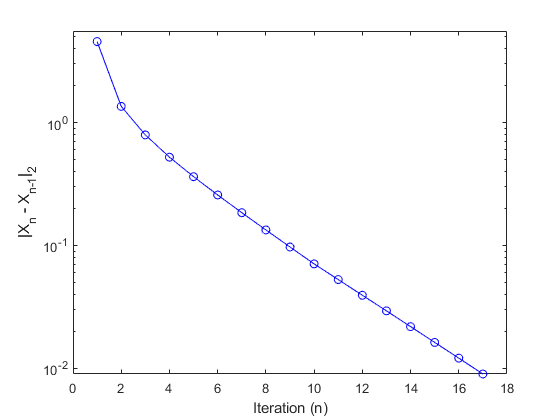}
    \caption{The convergence of the initial mesh smoothing.}
    \label{fig:conv}
  \end{subfigure}

  \begin{subfigure}{.5\textwidth}
    \includegraphics[width=.98\linewidth]{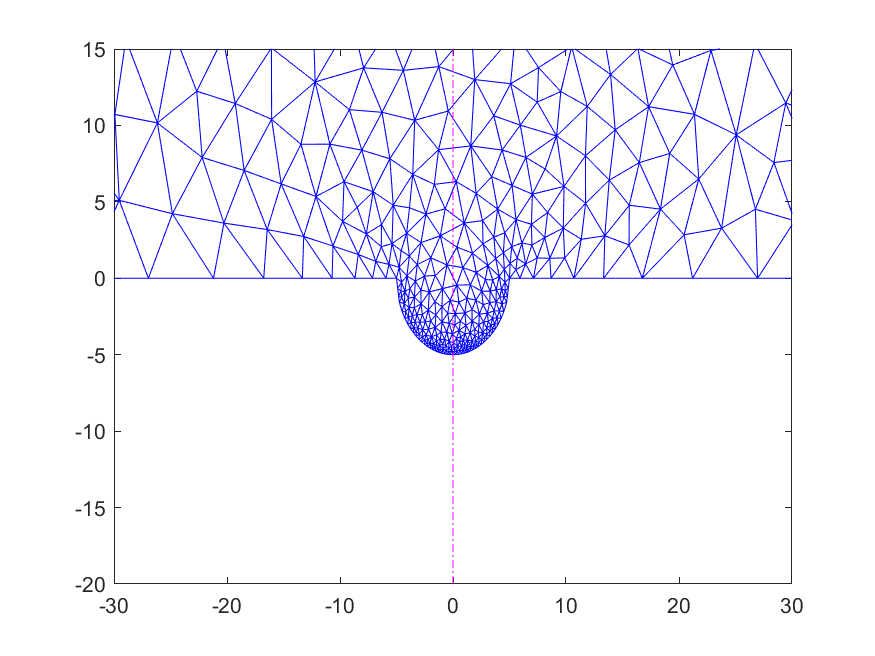}
    \caption{The smoothed initial mesh.}
    \label{fig:smooth}
  \end{subfigure}
  \begin{subfigure}{.5\textwidth}
    \includegraphics[width=.98\linewidth]{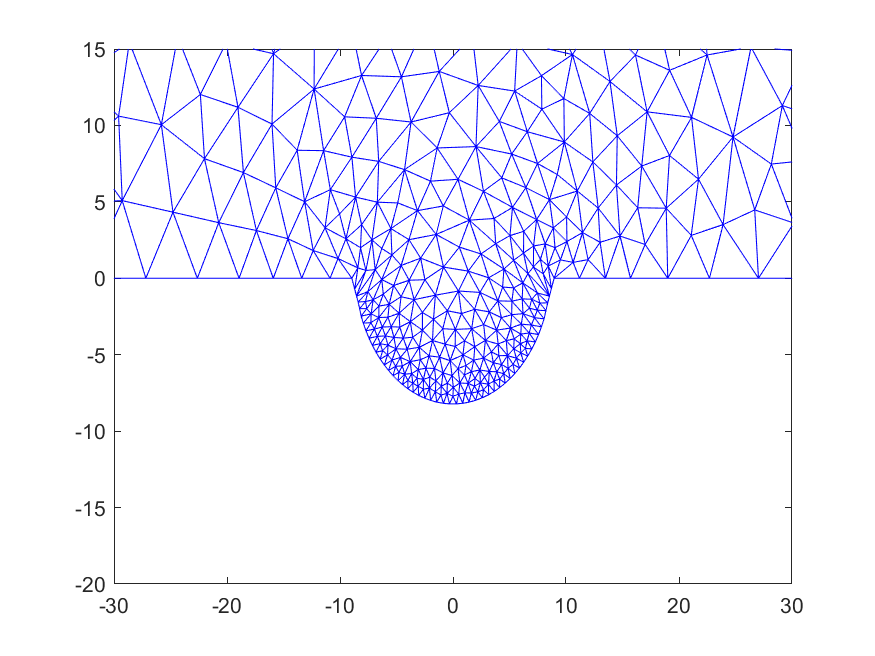}
    \caption{The mesh after $t=60$ s.}
    \label{fig:finalmesh}
  \end{subfigure}
  \caption{(a) A nonuniform initial mesh, (b) the effect of mesh smoothing on the positions of the nodes, (c) the mesh after smoothing,  and  (d) the mesh a $t=60$ s  using the monitor function \eqref{eq_macKenzie_mf_2}.}
  \label{fig:non-uniform-mesh-over-time}
\end{figure}

Presented in Figure \ref{fig:nonuniform} is the non-uniform initial mesh generated with the
Matlab function {\bf initmesh} where 45 nodes are located on the pit boundary. In order to optimize this mesh, smoothing steps were performed and the sum of the absolute differences in position of the nodes between subsequent smoothing iterations is displayed in Figure \ref{fig:conv}. The motion of the nodes decreases with iteration number and after 17 iterations the absolute difference is down to $10^{-2}$. The optimized mesh is shown in Figure \ref{fig:smooth} and significant differences in the locations of the nodes both inside and outside the pit are observed. Note that this smoothing operation only needs to be performed once since the results can be saved and used as the starting mesh in subsequent experiments. The computational mesh after the pit has evolved for 60 s is presented in Figure \ref{fig:finalmesh} and it is observed that node spacing within the pit remains very good.

The results in Figures \ref{fig:uniform-mesh-over-time} and  \ref{fig:non-uniform-mesh-over-time} show
that the MMPDE approach is robust with respect to the initial grid, continuously evolving the mesh according to changing domain and solution features.
Even with a uniform initial mesh, the MMPDE approach does quite well,
automatically recovering the requested increased mesh density near the pit
boundary.  We do notice, however, some additional stretching of the nodes in
this case as compared to the simulation which starts from an improved
non-uniform initial mesh.   The stretching of the mesh can be reduced through the use of monitor functions designed to control the shape of the elements.

The relatively small scale mesh smoothing is particularly useful should the pit boundary movement be large or if the pit boundary movement induces a discontinuous change in the geometry (during pit merging for example).  We will see this in Section \ref{sec:numerical_results}.

\subsection{Effect Of $\tau$ On The Moving Mesh}

\begin{figure}[ht!]
  \begin{subfigure}{.33\textwidth}
    \centering
    \includegraphics[width=.98\linewidth,height=.65\linewidth]{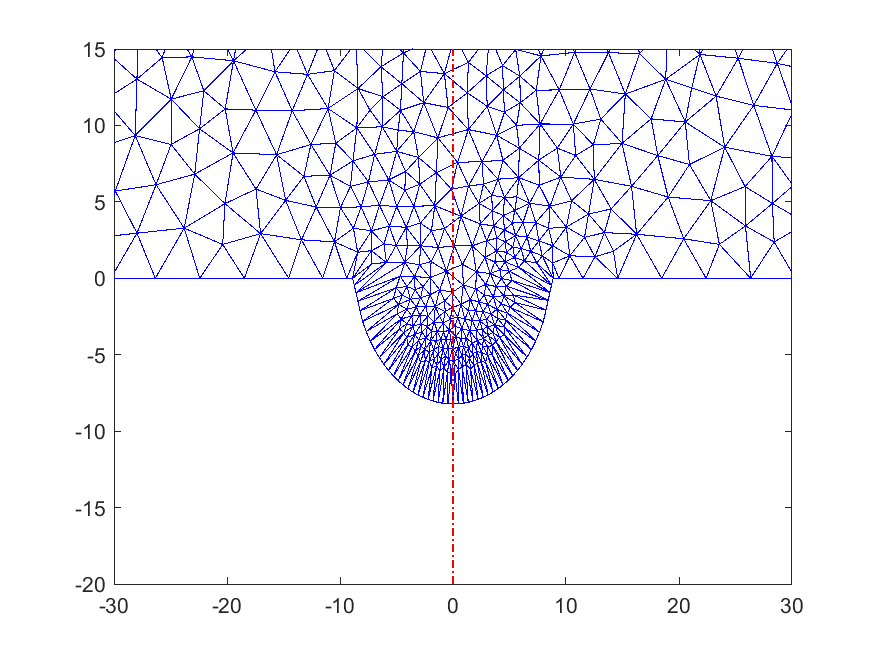}
    \caption{Using $\tau = 10^{-2}$.}
  \end{subfigure}
  \begin{subfigure}{.33\textwidth}
    \centering
    \includegraphics[width=.98\linewidth,height=.65\linewidth]{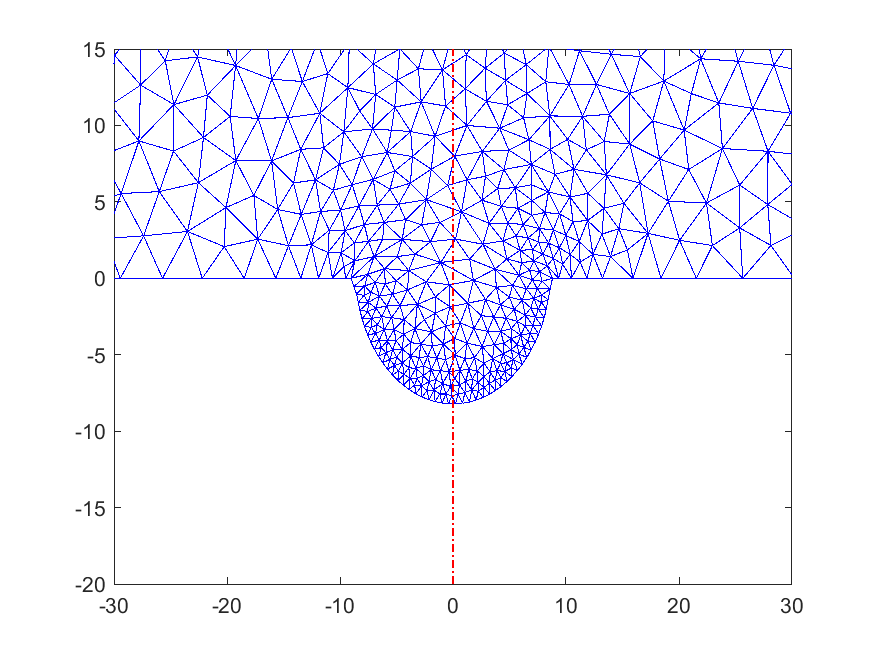}
    \caption{Using $\tau = 10^{-4}$.}
  \end{subfigure}
  \begin{subfigure}{.33\textwidth}
    \centering
    \includegraphics[width=.98\linewidth,height=.65\linewidth]{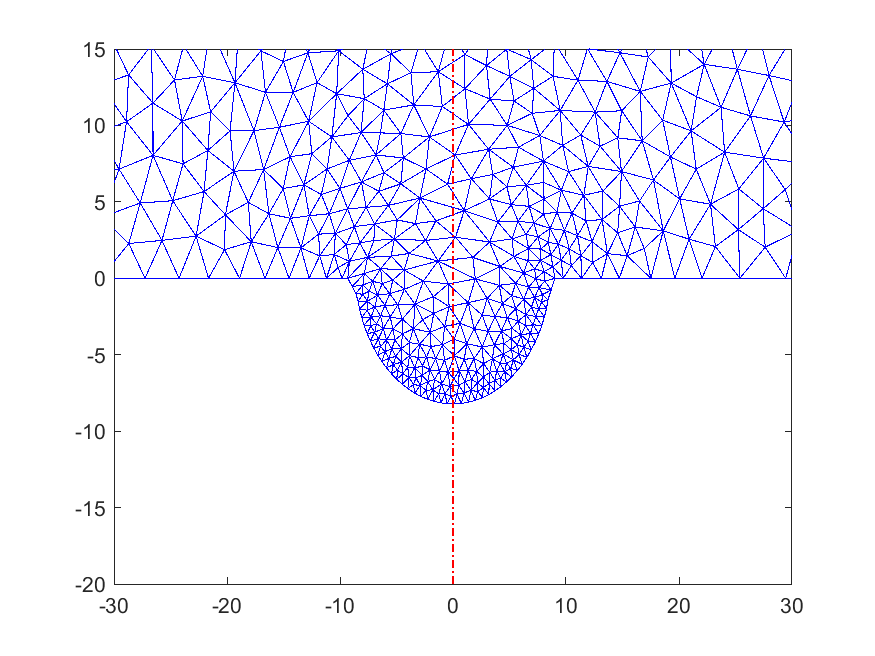}
    \caption{Using $\tau = 10^{-6}$.}
  \end{subfigure}
  \caption{Effect of $\tau$ on the mesh after 120 s with the monitor function
    \eqref{eq_macKenzie_mf_2} using $\mu_1 = 100$ and $\mu_2=1$.}
  \label{fig:effect-of-tau}
\end{figure}

At the end of Section 2.3 we mentioned the MMPDE (relaxation) parameter $\tau$ and here we
demonstrate the effect of $\tau$ on the moving mesh. 
Instead of forcing exact equidistribution each time $t$, we relax
the condition and require equidistribution at time $t+\tau$.  Hence, the smaller
the size of $\tau$ the quicker the mesh will react to changing features in the
solution.  To demonstrate this effect, we start from a uniform initial grid and
show the resulting meshes after $t=60$ s using three values of $\tau$, as shown in Figure~\ref{fig:effect-of-tau}.  We observe a greater concentration of nodes near the pit boundary for smaller values of $\tau$. Larger values of $\tau$ lead to stretched elements; the mesh is not able to keep up with the changing computational domain. This relaxation does come at a cost, however, as smaller values of $\tau$ require more time steps for the integration of the MMPDE.  In practice, one should select $\tau$ in tandem with mesh density parameters, choosing the largest value of $\tau$ which allows a balance of computational cost and mesh quality.

\subsection{Alternating Mesh and Physical PDE Iteration}

\begin{figure}[ht]
  \centering
  \begin{subfigure}[b]{0.7\textwidth}
  \centering
  \tikzset{>={Latex[width=3mm,length=3mm]}}  
  \begin{tikzpicture}
    \begin{scope} [scale = 0.28mm]
      \def\w{5}      

      \draw[->,gray, very thick] (-1.9*\w/2,\w/4) -- (-\w/2,\w/4);
      \draw[->,gray, very thick] (\w/2,\w/4) -- (1.9*\w/2,\w/4);
      \draw[dashed,gray, very thick] (-\w/2,\w/4) -- (\w/2,\w/4);
      \draw (-\w/2,0) -- (\w/2,0) -- (\w/2,\w/2) -- (-\w/2,\w/2) -- (-\w/2,0.0) -- (-\w/2,0.0);

      \draw[blue] (-1.45*\w/2,\w/4+0.85) node [align=center, below]{\footnotesize{ $\mathbf{x}^n$, $\mathbf{\varphi}^n$}};
      \draw[blue] (-0.05,\w/4+0.7) node [align=center, below]{\footnotesize{ Adaptive Mesh Generator}};
      \draw[blue] (0,\w/4-0.2) node [align=center, below]{\footnotesize{ PDE Solver}};
      \draw[blue] (1.45*\w/2,\w/4+0.85) node [align=center, below]{\footnotesize{ $\mathbf{x}^{n+1}$, $\mathbf{\varphi}^{n+1}$}};

    \end{scope}
  \end{tikzpicture}
  \caption{\label{fig:sim} The simultaneous solution approach.}
\end{subfigure}
\begin{subfigure}[b]{0.7\textwidth}
  \centering
  \includegraphics[width=13.25cm]{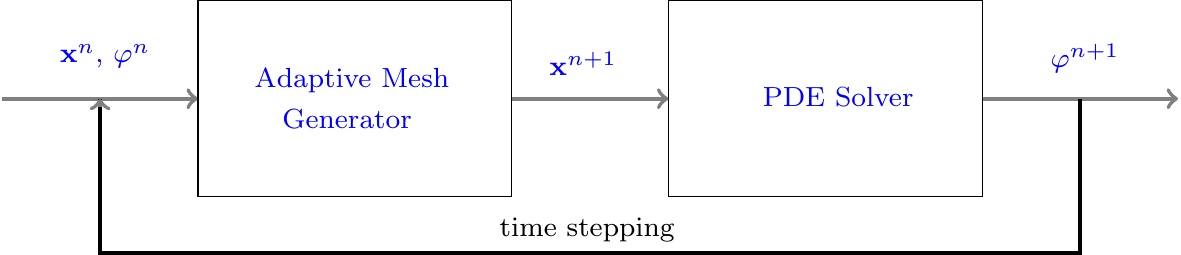}
  \caption{\label{fig:alt} The alternating solution approach.}
\end{subfigure}
\caption{The (a) simultaneous and (b) alternating approaches to solve the coupled corrosion model and mesh PDE.}
\end{figure}

There are two approaches that can be used to solve the coupled physical PDE and
mesh equation: a simultaneous or an alternating approach. In a simultaneous
approach, the discrete physical PDE and the discrete mesh equation provide a
fully coupled system for both the mesh and solution unknowns, as shown in
Figure~\ref{fig:sim}. The disadvantage of this approach is the highly nonlinear
coupling between the physical solution and the mesh, resulting in a potentially difficult, large discrete system.

On the other hand, the alternating solution approach generates the mesh
$\mathbf{x}^{n+1}$ at a new time step  using the physical solution $\varphi^n$
and the mesh $\mathbf{x}^n$ at the current time. Then   the solution
$\varphi^{n+1}$ at the new time level is computed,  as shown in
Figure~\ref{fig:alt}. In this approach, there may be a lag between the solution
and the mesh.  Generally, this does not create any difficulties if the time step
is reasonably small. The main advantages of the alternating approach are: (i)
the mesh generation code is not directly coupled to the physical PDE solve
thereby  increasing flexibility and reusability of code, (ii) the mesh PDE and physical PDE solvers can be developed and optimized in a modular way,
and hence (iii) the individual mesh and physical PDE solvers are more efficient. 
MMPDElab uses this alternating approach.

\subsection{Solution Of The Moving Boundary Value Problem}
\begin{figure}[ht!]
  \centering
  \includegraphics[height=2.20cm]{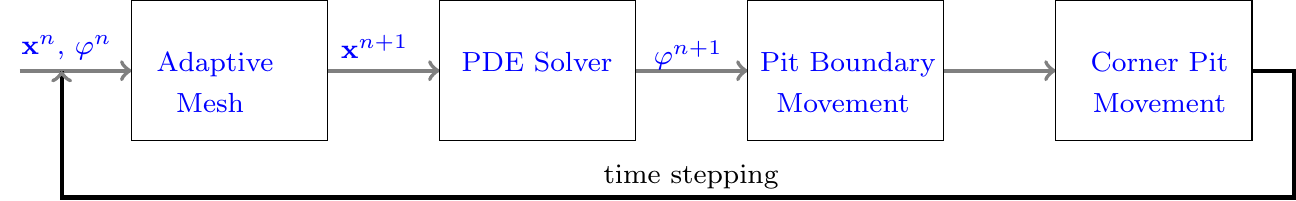}
  \caption{Flow chart for the physical PDE solve, the    mesh PDE solve, and the pit boundary movement.}
  \label{fig:flowchart}                                \end{figure}

The flowchart in Figure \ref{fig:flowchart} outlines the implementation of our
computational pitting corrosion model using MMPDElab framework. The first two
steps are the same as the alternating step approach given in Figure
\ref{fig:alt}. The pit boundary is then moved based on the new positions of the
adaptive mesh, followed by the movement of the corners of the pit. These last steps are detailed in the following section.

\subsection{Details Of The Pit Boundary Movement} \label{sec:pit_movement}
To get the new position of the pit we have to specify a direction
and magnitude of movement for each node on the boundary of the pit and the appropriate movement for each corner of the pit.  A pit corner is a vertex which is part of the pit boundary and has a $y$--coordinate of zero.

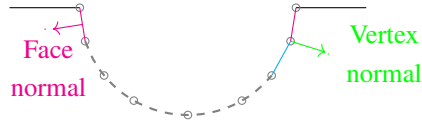
\begin{figure}[ht]
  \centering
  \begin{tikzpicture} [scale = 0.25mm]
    \draw (-3.3,0) -- (-2,0);
    \draw (2,0) -- (3.3,0);
    \draw [gray, dashed, thick, domain=-1.9:1.55] plot (\x, {-sqrt(4-\x*\x)});

    \draw[magenta] [->] (-1.95,-0.312) -- (-2.5,-0.4);    
    \filldraw[magenta]  (-2.64,-0.4) circle (0.1pt) node[align=center, below]{\footnotesize{ Face }\\ \footnotesize{ normal}};

    \foreach \x in {-2,-1.9, -1.55, -1, 0, 1., 1.55,1.9,2}  
    \draw[gray]  (\x, {-sqrt(4-\x*\x)}) circle (0.07);
    \draw[magenta]  (-2,0)-- (-1.9,{-sqrt(4-1.9*1.9)});
    \draw[magenta] (1.9,{-sqrt(4-1.9*1.9)}) -- (2,0);
    \draw[cyan] (1.55,{-sqrt(4-1.55*1.55)}) -- (1.9,{-sqrt(4-1.9*1.9)});
    \draw[green] [->] (1.9,{-sqrt(4-1.9*1.9)}) -- (2.55,-0.83);   
    \filldraw[green]  (2.6,-0.86) circle (0.1pt) node[align=center, right]{\footnotesize{ Vertex}\\ \footnotesize{ normal}};
  \end{tikzpicture}
  \caption{Definitions of face and vertex normals.} \label{fig:fv_normals}
\end{figure}

As shown in Figure~\ref{fig:fv_normals}, a face normal is the outward pointing vector perpendicular to an edge or segment joining vertices.  Taking the average of two face normals on adjacent edges gives us the vertex normal for the vertex between those edges.
The vertex normals give the direction of movement for the pit boundary.

As mentioned in Section 2.1, the magnitude of the normal velocity of each vertex on the boundary of the pit is given by (\ref{eq:vn}), which we may write as
\begin{equation}\label{eq:normal_velocity}
  V_{n} = \frac{1}{ c_\text{solid}}\cdot A_\text{diss}\cdot e^{\Big( \frac{z F ( V_{\text{corr}} + \alpha( V_\text{app} - V_{\text{corr}} - \varphi))}{RT}  \Big)}.
\end{equation}
Once the vertices on the pit boundary are moved, the location  of the corner nodes for
the pit are updated using the following procedure.  A linear extrapolation of
the edge joining the two vertices that are closest to the corner and lie on the
pit boundary is computed.   The new corner location is given by the intersection
of this line and $y=0$, as necessary.

\begin{figure}[h]
  \begin{subfigure}{.5\textwidth}
    \centering
    \includegraphics[width=.9\textwidth]{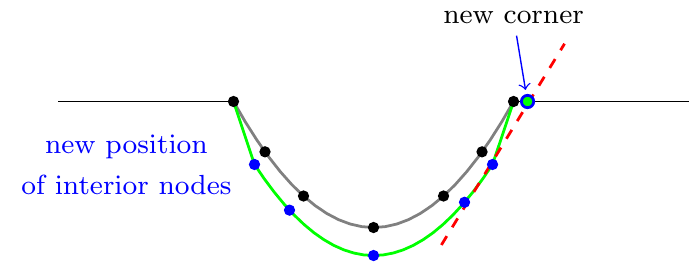}
    \caption{New corner remains on $y=0$.}
\label{fig:nomove}
  \end{subfigure}                                           \begin{subfigure}{.5\textwidth}                                   \centering

 \includegraphics[width=0.9\textwidth]{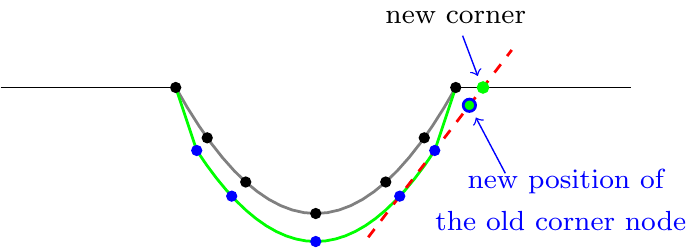}
    \caption{Corner node movement onto the pit boundary.}
    \label{fig:cornermove}
  \end{subfigure}
  \caption{Updating the corner position: (a) the corner is moved to its new location along $y=0$ or (b) the (old) corner is moved onto the boundary of the pit. }
  \label{fig:newcorner}
\end{figure}

There are two situations which may arise as shown in Figure \ref{fig:newcorner}.  If the new corner is close to the old corner (Figure \ref{fig:nomove}), then no further changes are required.   If the new corner is not close to the old corner (Figure \ref{fig:cornermove}) then the old corner is moved into the pit using the same extrapolation line. The idea here is to support large movements of the pit boundary by moving grid points along $y=0$ into the pit boundary.

\subsection{Merging Pits}\label{sec:merging}
\begin{figure}[ht]
  \centering  
  \includegraphics[width= 0.99\textwidth]{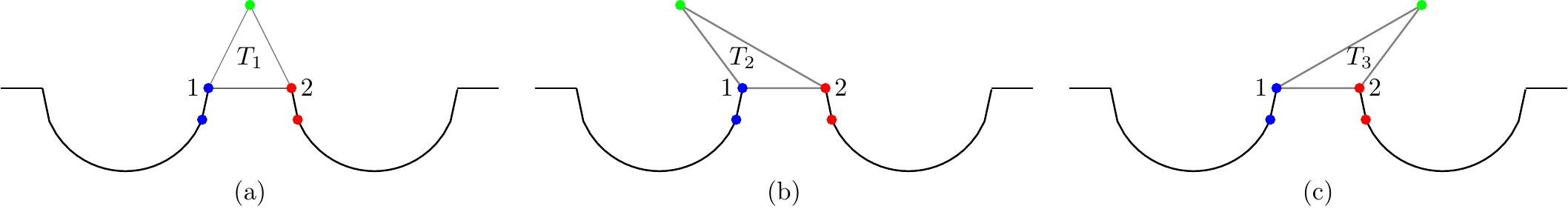}
  \caption{Three possible element orientations between the pits at the time of a
    merge}
  \label{fig:pit_merge2}
\end{figure}

\begin{figure}[ht!]
  \centering
  \begin{subfigure}{.30\textwidth}
    \centering
    \includegraphics[width=.98\linewidth]{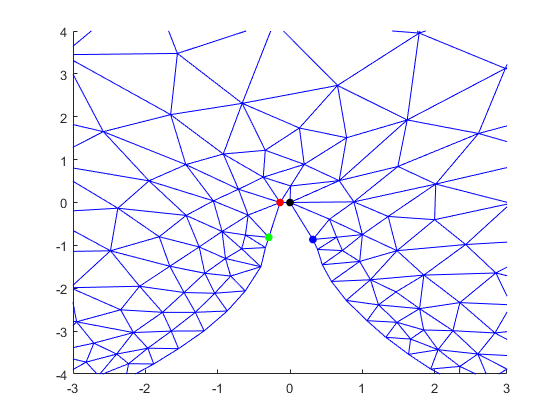}
    \caption{A pit merge is initiated.}
    \label{fig:pit-merge-a}
  \end{subfigure}
  \begin{subfigure}{.3\textwidth}
    \centering
    \includegraphics[width=.98\linewidth]{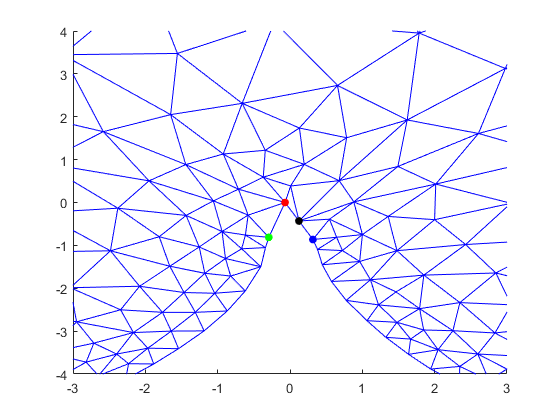}

    \caption{A merge with no mesh topology change.}
    \label{fig:pit-merge-b}
  \end{subfigure}
  \begin{subfigure}{.3\textwidth}
    \centering
    \includegraphics[width=.98\linewidth]{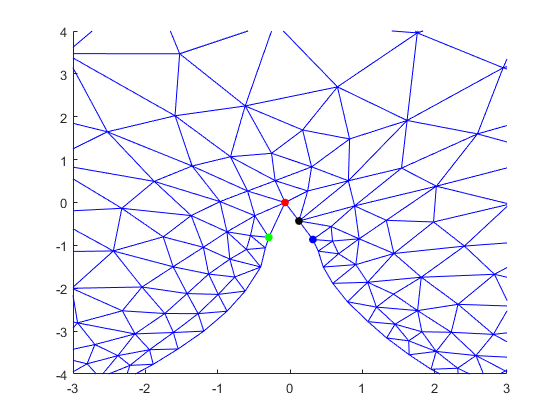}

    \caption{Post merge mesh smoothing.}
    \label{fig:pit-merge-c}
  \end{subfigure}
  \caption{The pit merging process.}
  \label{fig:pit-merge}
\end{figure}

As multiple pits evolve, two pits may merge and form a single, larger pit. 
We visualize the pit merging procedure in Figure~\ref{fig:pit-merge}.
The merging process is initiated (Figure \ref{fig:pit-merge-a}) when there is a
single edge between two pits that is less than the user prescribed tolerance.
The left and right endpoints of that edge are tagged with red and  black,
respectively. In Figure \ref{fig:pit-merge-b} these two points merge to a single
point.  In order to avoid changing the number of mesh nodes and mesh topology,
either the black or the red node in Figure \ref{fig:pit-merge-a} has to move
into the pit boundary either on the left or the right side of the apex.  For
example, the black and red nodes can merge, creating a new red node (see Figure
\ref{fig:pit-merge-b}),  and the black node will move half-way between the red
and green nodes.  Alternately, the black and red node merge, creating a new
black node and the red node moves half-way between the red node and green node
(not shown).  We choose the vertex corresponding to the larger angle in the
element whose bottom edge is the single edge between the pits, see Figure
\ref{fig:pit_merge2}.   Once the merge has occurred, a mesh smoothing procedure
is used to obtain the mesh shown in Figure \ref{fig:pit-merge-c}. We can see that the mesh smoothing has evened out the size and shape of the elements to the right of the red apex node.

For each time step after the merge the location of the
apex is obtained as the intersection point of the linear extrapolations of the second last
edges to the left and right of the apex. 

\section{Numerical Results}\label{sec:numerical_results}
Based on the simple experiments in Section 3, throughout this numerical results section a constant value of $\tau = 10^{-5}$ is used in the MMPDE, and the constants  $\mu_1 = 100$ and $\mu_2 = 1$  are used as  default values for the mesh density function \eqref{eq_macKenzie_mf_2}.
The initial number of mesh points inside the pit is set to $61$. 
\subsection{Single Pit Simulations }
We begin by using our computational framework to compare the evolution of a single pit in three cases:
a homogeneous solid material (without a crystal direction),  a solid material
with a specified crystal direction, and a solid material with a discontinuity in
the crystal direction.   All simulations use an initial mesh constructed by
smoothing the non-uniform mesh generated using {\bf initmesh} as discussed in Section 3.3.

\begin{figure}[ht!]
  \begin{subfigure}{.5\textwidth}
    \centering
    \includegraphics[width=.98\linewidth]{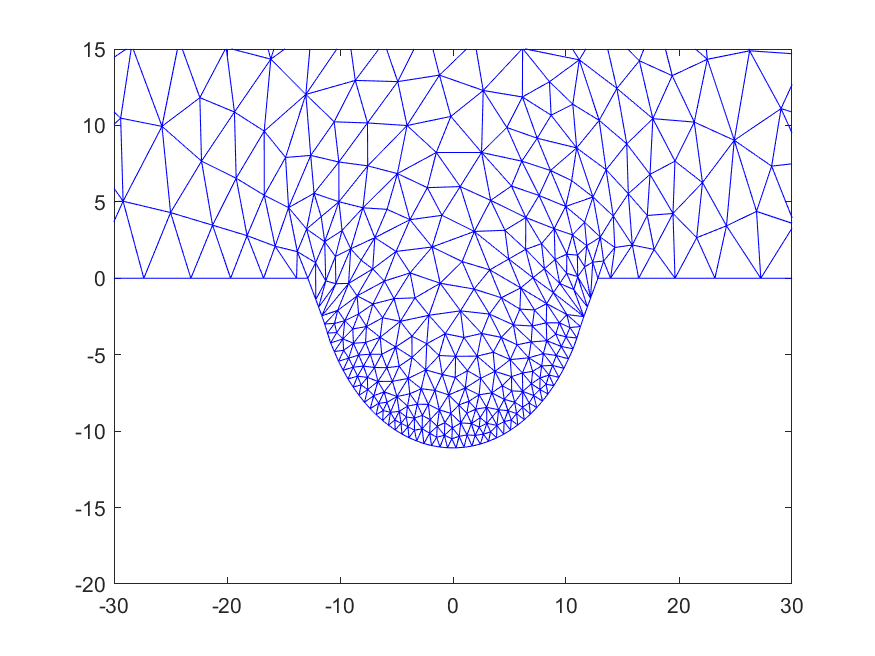}
    \caption{Mesh for a homogeneous crystal.}
    \label{fig:single-pit-b}
  \end{subfigure}
  \begin{subfigure}{.5\textwidth}
    \centering
    \includegraphics[width=.98\linewidth]{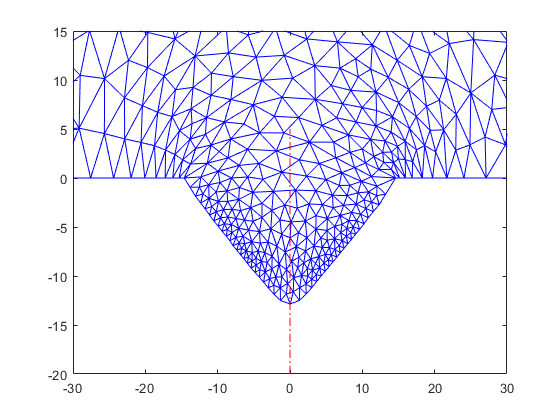}
    \caption{Mesh for a crystal with direction [001].}
    \label{fig:single-pit-e}
  \end{subfigure}
  \begin{subfigure}{.5\textwidth}
    \centering
    \includegraphics[width=.98\linewidth]{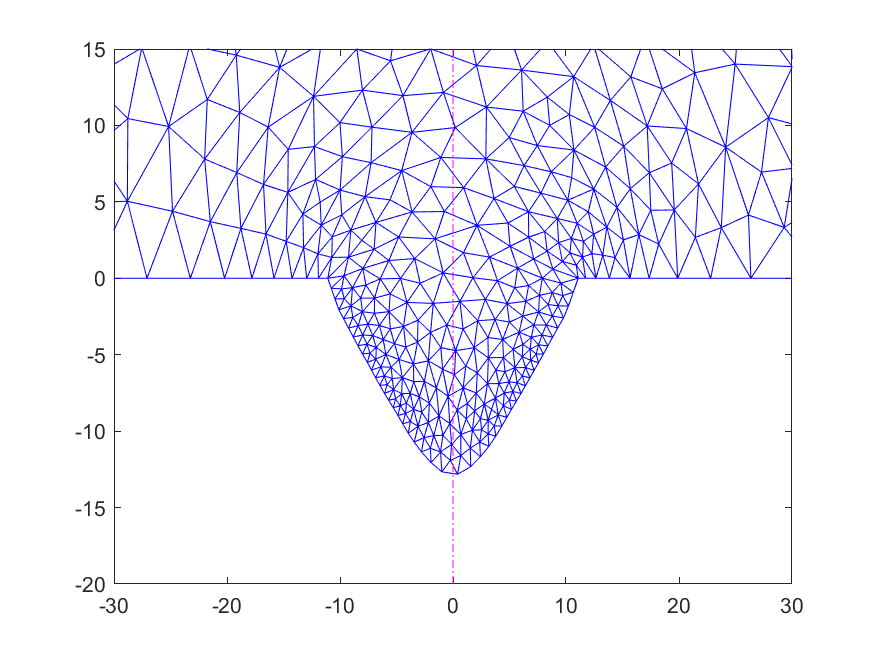}
    \caption{Mesh for a crystal with direction [101].}
    \label{fig:single-pit-c}
  \end{subfigure}
  \begin{subfigure}{.5\textwidth}
    \centering
    \includegraphics[width=.98\linewidth]{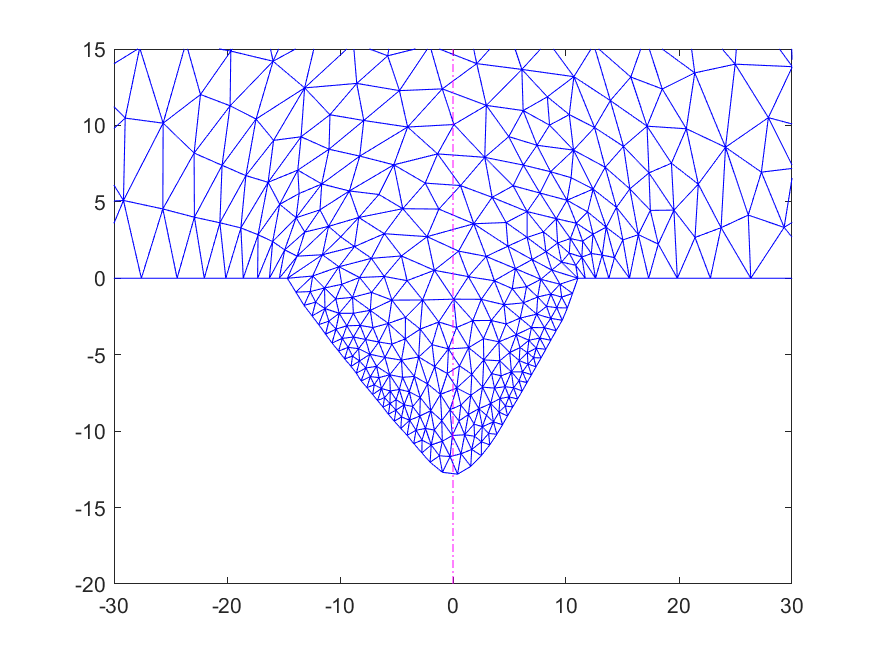}
    \caption{Mesh for a crystal with two directions [001] and [101].}
    \label{fig:single-pit-d}
  \end{subfigure}
  \caption{Pit configurations and meshes at $t=120$ s  for a) a homogeneous
    material, b) a single crystal oriented with a zone axis along [001], c) a
    single crystal oriented with a zone axis along [101], and d) a crystal with
    an interface at $x=0$; the crystal directions to the left and right of $x=0$ are [001] and [101], respectively.}
  \label{fig:single-pit}
\end{figure}

Figures~\ref{fig:single-pit-b} 
\ref{fig:single-pit-e},
\ref{fig:single-pit-c} and \ref{fig:single-pit-d}   show the
final pit geometry and final meshes for a homogeneous material, a material with crystal direction $[0 0
1]$, a material with  crystal direction $[1 0 1]$,  and a material with a
discontinuity in crystal directions, $[0 0 1]$ to the left of $x=0$ and $[1 0
1]$ to the right of $x=0$.


It is observed in Figure \ref{fig:single-pit-b} that the final shape of the
homogeneous pit is the same as the initial pit since the chosen $V_{corr}$ is a
constant value of 
$-0.24\ \textrm{V}$, that is, the same in every direction within the pit. Hence, from equation (\ref{eq:normal_velocity}) the normal velocity at all locations within the pit will be equal. The situation is not the same for pits with a crystalline structure since $V_{corr}$ will vary with crystal direction. For example, for crystal directions of the forms $\left\langle 001\right\rangle $, $\left\langle 011\right\rangle $ and $\left\langle 111\right\rangle $, $V_{corr}$ will have values of -0.2297 V, -0.2455 V, and -0.2525 V, respectively. Again from equation (\ref{eq:normal_velocity}), we see that the normal velocity is greater for lower magnitude $V_{corr}$ values; that is, $V_n\left( 111\right) < V_n\left( 011\right)  < V_n\left( 001\right) $. The effect of this $V_{corr}$ dependent velocity is displayed in Figure \ref{fig:single-pit-e} where we observe that the sides of the pit have become straight and angled 90 degrees with one another (when axes are equally scaled). This behaviour is expected. As shown in Figure \ref{fig:fig2}, for a crystal oriented with a zone axis along $\left[ 001\right] $, $\left\langle 001\right\rangle$ directions are located along the horizontal and vertical axes and $\left\langle 011\right\rangle$ directions midway in between. Thus, we expect that the $\left[ 100\right]$ and $[0\bar{1}0]$ directions will move faster than the $\left[ 1\bar{1}0\right]$ direction. As the faster locations on the pit boundary move, their orientation will change and eventually become the same direction as the slowest moving axis, in this case $\left[ 1\bar{1}0\right]$. For all future times, the sides of the pit will move outward perpendicular to these two lines while maintaining the same angular relationship. We observe the same effect in Figure \ref{fig:single-pit-c} where the crystal has been oriented along a zone axis of $\left[ 101\right]$. In this case, the slowest moving directions are along  $\left\langle 111\right\rangle$ and the angle between the $\left[ \bar{1}\bar{1}1\right]$ and $\left[ 1\bar{1}1\right]$ planes agrees with the expected value of 70.5 degrees. Figure \ref{fig:single-pit-d} displays the final pit shape where there is a discontinuity in the crystal directions and the left and right sides of the crystal were oriented along zone axes of $\left[ 001\right]$ and $\left[ 101\right]$, respectively. The left and right sides of the pit are straight lines moving along $\left[ \bar{1}\bar{1}0\right]$ and $\left[ 1\bar{1}1\right]$ directions, respectively.

\subsection{Multiple Pit Simulations }
As mentioned previously, if multiple pits exist in a material and the pits grow
large enough, there is the potential that the pits will merge during the
simulation.
The  imminent merge needs to be detected, the boundaries of each of the previously isolated pits need to be updated in a way which avoids boundary crossing, and the mesh around merge location needs to be adjusted smoothly and without topology changes.
Once complete the merged pit is then treated as one larger pit, and the
evolution continues. See Section \ref{sec:merging} for details.

\begin{figure}[ht!]
  \begin{subfigure}{.5\textwidth}
    \centering
    \includegraphics[width=.98\linewidth]{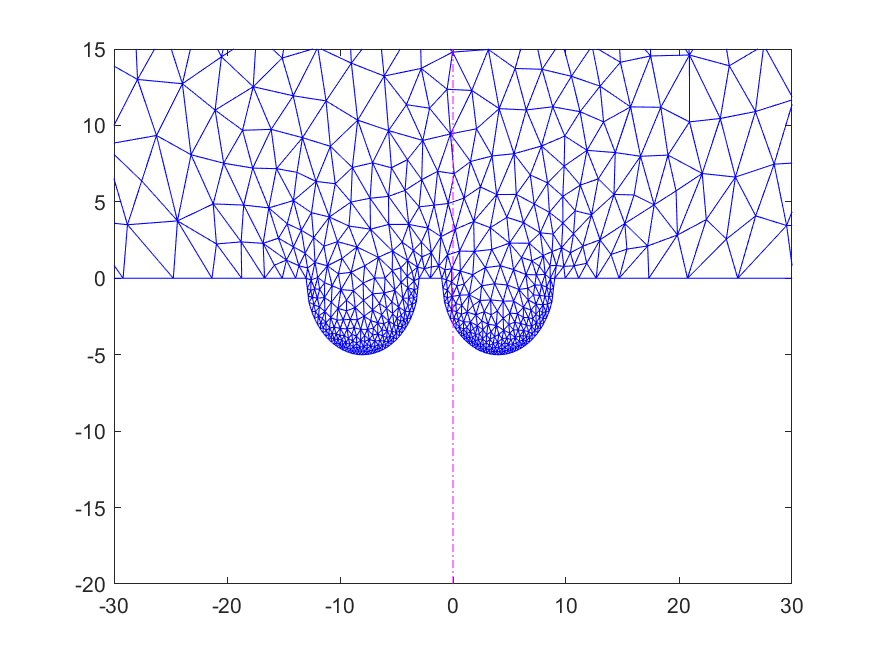}
    \caption{Initial mesh for two pits. }
    \label{fig:multiple-pit-a}
  \end{subfigure}
  \begin{subfigure}{.5\textwidth}
    \centering
    \includegraphics[width=.98\linewidth]{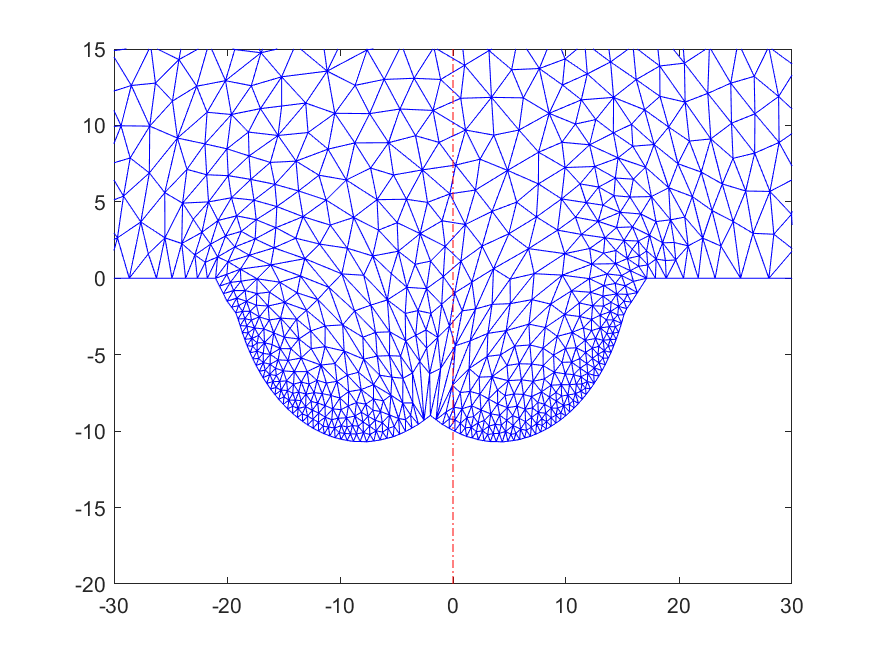}
    \caption{Mesh for merged pits at $t=120$ s.}
    \label{fig:multiple-pit-b}
  \end{subfigure}
  \begin{subfigure}{.5\textwidth}
    \centering
    \includegraphics[width=.98\linewidth]{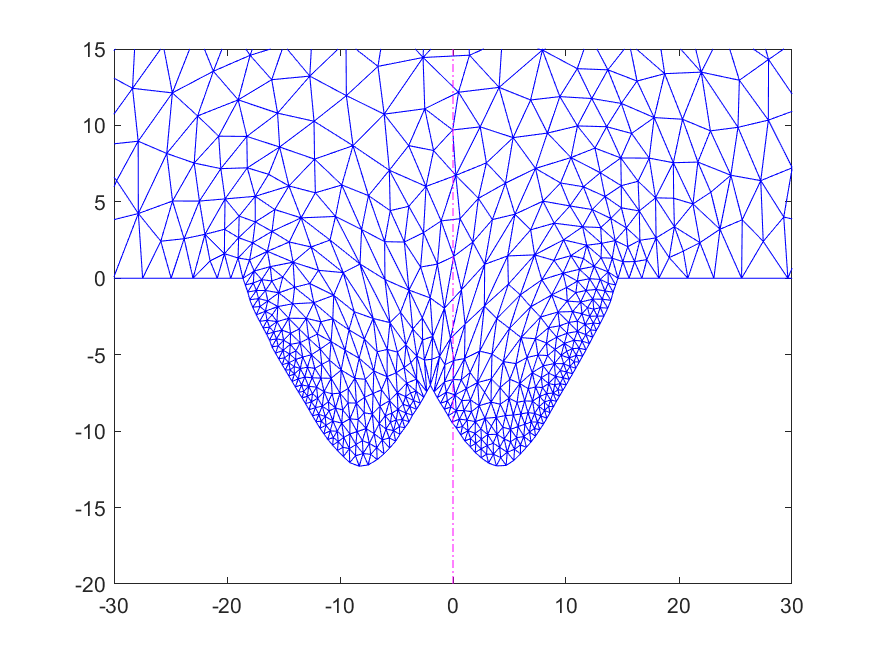}
    \caption{Mesh for merged pits for a material with a single crystal
      direction $[1 0 1]$ at $t=120$ s.}
    \label{fig:multiple-pit-c}
  \end{subfigure}
  \begin{subfigure}{.5\textwidth}
    \centering
    \includegraphics[width=.98\linewidth]{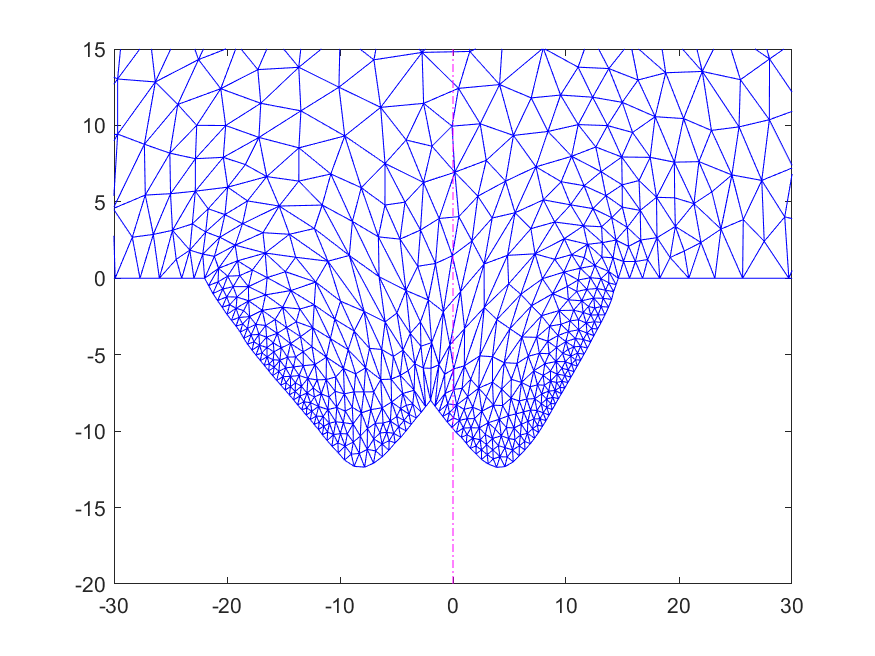}
    \caption{Mesh for merged pits for a material with two crystal
      directions, $[0 0 1]$ if $x < 0$ and $[1 0 1]$ if $x > 0$, at $t=120$ s. }
    \label{fig:multiple-pit-d}
  \end{subfigure}
  \caption{Pit evolution and adaptive mesh generation for merging multiple pits for three material configurations. }
  \label{fig:multiple-pit}
\end{figure}

To demonstrate the robustness of our adaptive simulation framework for the
evolution of multiple pits we start with two pits relatively close together in a
homogeneous material, as shown in plot Figure~\ref{fig:multiple-pit-a}.  A
quality initial mesh which concentrates nodes near the boundary of both pits is
generated as discussed in Section \ref{sec:initialmeshgeneration}.  As the corrosion continues the pits grow, and
hence grow closer together.   The pits then merge and continue to evolve as
shown in Figure~\ref{fig:multiple-pit-b}.
Figure~\ref{fig:multiple-pit-c} provides the result of a similar simulation with
a material oriented in the $[1 0 1]$ crystal direction.    The crystal direction
clearly affects the geometry of the merged pit.  Figure~\ref{fig:multiple-pit-d}
shows the resulting pit geometry and associated mesh for  merged pits in a material with two crystal directions, where the discontinuity in crystal direction is located at $x=0$.


  \begin{figure}

  \begin{subfigure}{.5\textwidth}
    \centering
    \includegraphics[width=.98\linewidth]{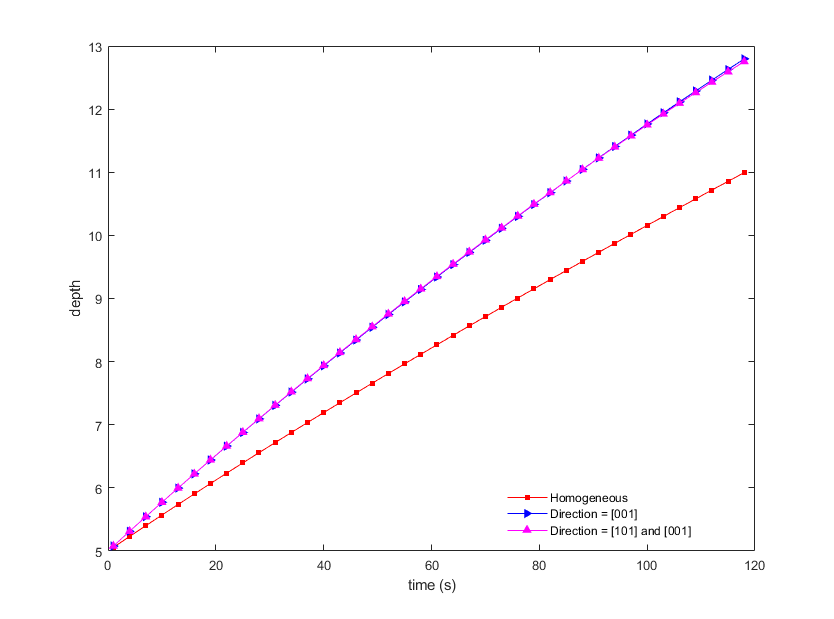}
    \caption{Pit-depth over time for homogeneous  and non-homogeneous crystals.}
    \label{fig:depthhomononhomo}
  \end{subfigure}
  \begin{subfigure}{.5\textwidth}
    \centering
    \includegraphics[width=.98\linewidth]{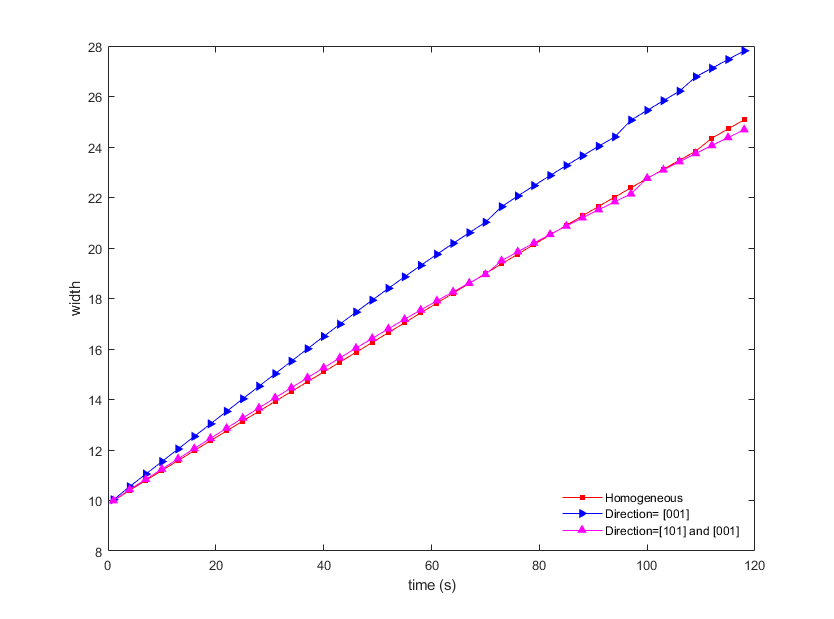}
    \caption{Pit-width over time for homogeneous and non-homogeneous crystals.}
    \label{fig:widthhomononhomo}
  \end{subfigure}

  \caption{Pit-depths and widths for homogeneous and non-homogeneous crystals.}
  \label{fig:???}
\end{figure}

We now more closely study the pit depth and width as a function of time for a
single homogeneous pit, a single crystal oriented with a zone axis along [001],
and two crystals with an interface at $x=0$ where the left and right crystals
are oriented along [001] and [101] zone axes, respectively.  Recall, the initial and
final pit
configurations for these three situations are displayed in Figures
\ref{fig:single-pit-b}, \ref{fig:single-pit-e} and \ref{fig:single-pit-d},
respectively.  Plotting pit depth and width as a function of time leads to
nonlinear curves as shown in Figures \ref{fig:depthhomononhomo} and
\ref{fig:widthhomononhomo}. It has been common practice to fit corrosion loss
curves using a power-law equation and for the initial stages of corrosion it
seems to work well, see \cite{Melchers2019-pb}.
Since the loss of material due to corrosion is a function of the dimensions of the pit, it is expected that the same power-law behaviour should hold for our pit depth and width data.  The model we use is

\begin{equation*}
	\textrm{depth(t)}\ \left( \textrm{or}\right) \ \textrm{width(t)} = at^b +c,
\end{equation*}
where $a$, $b$ and $c$ are fitting parameters and the initial pit is defined when $t=1$. 
The curve fits were excellent and the fitting parameters found are presented in Table 2.  The initial pit width and depth were 10 microns and 5 microns, respectively, and these values are close to the value of $a+c$; the initial dimension of the pit predicted from the fitting procedure.  It is reassuring to note that the modelled pitting corrosion behaviour follows an expression used to fit experimental corrosion losses.

\begin{table}[h!]
	\begin{center}		
		\label{tab:fits}
		\begin{tabular}{lccc|ccc}
			\toprule
			\multirow{2}{*}{\textbf{Case}} &
                                       \multicolumn{3}{c}{\textbf{Width}} &
                                                                            \multicolumn{3}{c}{\textbf{Depth}} \\
                                     & a & b & c & a & b & c \\
      \hline
      homogeneous   & 0.142(2) & 0.980(3) & 9.83(2) & 0.076(1) & 0.917(3) & 4.95(1)\\
      {[001]}       & 0.243(8) & 0.907(6) & 9.62(5) & 0.116(3) & 0.886(4) & 4.89(2)\\
      {[001]/[101]} & 0.181(5) & 0.927(6) & 9.74(4) & 0.121(3) & 0.877(6) & 4.89(2)\\
			\bottomrule
		\end{tabular}
		\caption{Power-law model fitting parameters for the 6 curves presented in Figure 16. The numbers in brackets represent uncertainty in the last significant digit.}  
	\end{center}
\end{table}

\section{Conclusion}
We have presented a robust, fully automatic, moving mesh solution framework for pitting corrosion. The moving mesh approach continuously and smoothly evolves a fixed mesh topology according to changing pit geometry.  Single and multiple pits are considered, as are materials with different crystal direction(s).  A procedure is presented which allows pits to merge without a change in mesh topology, allowing computation to proceed without restarting the computation.

The simulation of large pit growth or the initiation of many pits would likely
benefit from an $hr$--refinement strategy (which both redistributes nodes as we
have presented here but also allows periodic changes to the number of mesh
nodes) coupled with a domain decomposition approach to allow the problem to be
spatially partitioned and the computation distributed to harness additional processors.  This will be the subject of future work.
Current work includes extending the computational framework to allow for
 more heterogenous materials,  with corrosive resistant ``pockets'' or holes
or voids.

\section{Data Availability}
The raw or processed data required to reproduce these figures and findings
cannot be shared at this time due to technical or time limitations, but are
available from the authors upon request.

\bibliographystyle{model1-num-names}
\bibliography{bibliography}
\end{document}